\catcode`\@=11
 
\font\teneufm=eufm10
\font\seveneufm=eufm7
\font\fiveeufm=eufm5
\newfam\eufmfam
\textfont\eufmfam=\teneufm  \scriptfont\eufmfam=\seveneufm
  \scriptscriptfont\eufmfam=\fiveeufm
\def\eufm@{\hexnumber@\eufmfam}

\font\tenmsa=msam10
\font\sevenmsa=msam7
\font\fivemsa=msam5
\font\tenmsb=msbm10
\font\sevenmsb=msbm7
\font\fivemsb=msbm5
\newfam\msafam
\newfam\msbfam
\textfont\msafam=\tenmsa  \scriptfont\msafam=\sevenmsa
  \scriptscriptfont\msafam=\fivemsa
\textfont\msbfam=\tenmsb  \scriptfont\msbfam=\sevenmsb
  \scriptscriptfont\msbfam=\fivemsb

\def\hexnumber@#1{\ifnum#1<10 \number#1\else
 \ifnum#1=10 A\else\ifnum#1=11 B\else\ifnum#1=12 C\else
 \ifnum#1=13 D\else\ifnum#1=14 E\else\ifnum#1=15 F\fi\fi\fi\fi\fi\fi\fi}

\def\msa@{\hexnumber@\msafam}
\def\msb@{\hexnumber@\msbfam}

\mathchardef\gx="2\eufm@78
\mathchardef\gg="2\eufm@67
\mathchardef\gm="2\eufm@6D
\mathchardef\gn="2\eufm@6E
\mathchardef\gd="2\eufm@64
\mathchardef\ga="2\eufm@61
\mathchardef\gh="2\eufm@68

\mathchardef\boxdot="2\msa@00
\mathchardef\boxplus="2\msa@01
\mathchardef\boxtimes="2\msa@02
\mathchardef\square="0\msa@03
\mathchardef\blacksquare="0\msa@04
\mathchardef\centerdot="2\msa@05
\mathchardef\lozenge="0\msa@06
\mathchardef\blacklozenge="0\msa@07
\mathchardef\circlearrowright="3\msa@08
\mathchardef\circlearrowleft="3\msa@09
\mathchardef\rightleftharpoons="3\msa@0A
\mathchardef\leftrightharpoons="3\msa@0B
\mathchardef\boxminus="2\msa@0C
\mathchardef\Vdash="3\msa@0D
\mathchardef\Vvdash="3\msa@0E
\mathchardef\vDash="3\msa@0F
\mathchardef\twoheadrightarrow="3\msa@10
\mathchardef\twoheadleftarrow="3\msa@11
\mathchardef\leftleftarrows="3\msa@12
\mathchardef\rightrightarrows="3\msa@13
\mathchardef\upuparrows="3\msa@14
\mathchardef\downdownarrows="3\msa@15
\mathchardef\upharpoonright="3\msa@16

\mathchardef\downharpoonright="3\msa@17
\mathchardef\upharpoonleft="3\msa@18
\mathchardef\downharpoonleft="3\msa@19
\mathchardef\rightarrowtail="3\msa@1A
\mathchardef\leftarrowtail="3\msa@1B
\mathchardef\leftrightarrows="3\msa@1C
\mathchardef\rightleftarrows="3\msa@1D
\mathchardef\Lsh="3\msa@1E
\mathchardef\Rsh="3\msa@1F
\mathchardef\rightsquigarrow="3\msa@20
\mathchardef\leftrightsquigarrow="3\msa@21
\mathchardef\looparrowleft="3\msa@22
\mathchardef\looparrowright="3\msa@23
\mathchardef\circeq="3\msa@24
\mathchardef\succsim="3\msa@25
\mathchardef\gtrsim="3\msa@26
\mathchardef\gtrapprox="3\msa@27
\mathchardef\multimap="3\msa@28
\mathchardef\therefore="3\msa@29
\mathchardef\because="3\msa@2A
\mathchardef\doteqdot="3\msa@2B

\mathchardef\triangleq="3\msa@2C
\mathchardef\precsim="3\msa@2D
\mathchardef\lesssim="3\msa@2E
\mathchardef\lessapprox="3\msa@2F
\mathchardef\eqslantless="3\msa@30
\mathchardef\eqslantgtr="3\msa@31
\mathchardef\curlyeqprec="3\msa@32
\mathchardef\curlyeqsucc="3\msa@33
\mathchardef\preccurlyeq="3\msa@34
\mathchardef\leqq="3\msa@35
\mathchardef\leqslant="3\msa@36
\mathchardef\lessgtr="3\msa@37
\mathchardef\backprime="0\msa@38
\mathchardef\risingdotseq="3\msa@3A
\mathchardef\fallingdotseq="3\msa@3B
\mathchardef\succcurlyeq="3\msa@3C
\mathchardef\geqq="3\msa@3D
\mathchardef\geqslant="3\msa@3E
\mathchardef\gtrless="3\msa@3F
\mathchardef\sqsubset="3\msa@40
\mathchardef\sqsupset="3\msa@41
\mathchardef\trianglerighteq="3\msa@44
\mathchardef\trianglelefteq="3\msa@45
\mathchardef\bigstar="0\msa@46
\mathchardef\between="3\msa@47
\mathchardef\blacktriangledown="0\msa@48
\mathchardef\blacktriangleright="3\msa@49
\mathchardef\blacktriangleleft="3\msa@4A
\mathchardef\blacktriangle="0\msa@4E
\mathchardef\triangledown="0\msa@4F
\mathchardef\eqcirc="3\msa@50
\mathchardef\lesseqgtr="3\msa@51
\mathchardef\gtreqless="3\msa@52
\mathchardef\lesseqqgtr="3\msa@53
\mathchardef\gtreqqless="3\msa@54
\mathchardef\Rrightarrow="3\msa@56
\mathchardef\Lleftarrow="3\msa@57
\mathchardef\veebar="2\msa@59
\mathchardef\barwedge="2\msa@5A
\mathchardef\doublebarwedge="2\msa@5B
\mathchardef\angle="0\msa@5C
\mathchardef\measuredangle="0\msa@5D
\mathchardef\sphericalangle="0\msa@5E
\mathchardef\varpropto="3\msa@5F
\mathchardef\smallsmile="3\msa@60
\mathchardef\smallfrown="3\msa@61
\mathchardef\Subset="3\msa@62
\mathchardef\Supset="3\msa@63
\mathchardef\Cup="2\msa@64

\mathchardef\Cap="2\msa@65

\mathchardef\curlywedge="2\msa@66
\mathchardef\curlyvee="2\msa@67
\mathchardef\leftthreetimes="2\msa@68
\mathchardef\rightthreetimes="2\msa@69
\mathchardef\subseteqq="3\msa@6A
\mathchardef\supseteqq="3\msa@6B
\mathchardef\bumpeq="3\msa@6C
\mathchardef\Bumpeq="3\msa@6D
\mathchardef\lll="3\msa@6E

\mathchardef\ggg="3\msa@6F

\mathchardef\circledS="0\msa@73
\mathchardef\pitchfork="3\msa@74
\mathchardef\dotplus="2\msa@75
\mathchardef\backsim="3\msa@76
\mathchardef\backsimeq="3\msa@77
\mathchardef\complement="0\msa@7B
\mathchardef\intercal="2\msa@7C
\mathchardef\circledcirc="2\msa@7D
\mathchardef\circledast="2\msa@7E
\mathchardef\circleddash="2\msa@7F
\def\ulcorner{\delimiter"4\msa@70\msa@70 }
\def\urcorner{\delimiter"5\msa@71\msa@71 }
\def\llcorner{\delimiter"4\msa@78\msa@78 }
\def\lrcorner{\delimiter"5\msa@79\msa@79 }
\def\yen{\mathhexbox\msa@55 }
\def\checkmark{\mathhexbox\msa@58 }
\def\circledR{\mathhexbox\msa@72 }
\def\maltese{\mathhexbox\msa@7A }
\mathchardef\lvertneqq="3\msb@00
\mathchardef\gvertneqq="3\msb@01
\mathchardef\nleq="3\msb@02
\mathchardef\ngeq="3\msb@03
\mathchardef\nless="3\msb@04
\mathchardef\ngtr="3\msb@05
\mathchardef\nprec="3\msb@06
\mathchardef\nsucc="3\msb@07
\mathchardef\lneqq="3\msb@08
\mathchardef\gneqq="3\msb@09
\mathchardef\nleqslant="3\msb@0A
\mathchardef\ngeqslant="3\msb@0B
\mathchardef\lneq="3\msb@0C
\mathchardef\gneq="3\msb@0D
\mathchardef\npreceq="3\msb@0E
\mathchardef\nsucceq="3\msb@0F
\mathchardef\precnsim="3\msb@10
\mathchardef\succnsim="3\msb@11
\mathchardef\lnsim="3\msb@12
\mathchardef\gnsim="3\msb@13
\mathchardef\nleqq="3\msb@14
\mathchardef\ngeqq="3\msb@15
\mathchardef\precneqq="3\msb@16
\mathchardef\succneqq="3\msb@17
\mathchardef\precnapprox="3\msb@18
\mathchardef\succnapprox="3\msb@19
\mathchardef\lnapprox="3\msb@1A
\mathchardef\gnapprox="3\msb@1B
\mathchardef\nsim="3\msb@1C
\mathchardef\napprox="3\msb@1D
\mathchardef\nsubseteqq="3\msb@22
\mathchardef\nsupseteqq="3\msb@23
\mathchardef\subsetneqq="3\msb@24
\mathchardef\supsetneqq="3\msb@25
\mathchardef\subsetneq="3\msb@28
\mathchardef\supsetneq="3\msb@29
\mathchardef\nsubseteq="3\msb@2A
\mathchardef\nsupseteq="3\msb@2B
\mathchardef\nparallel="3\msb@2C
\mathchardef\nmid="3\msb@2D
\mathchardef\nshortmid="3\msb@2E
\mathchardef\nshortparallel="3\msb@2F
\mathchardef\nvdash="3\msb@30
\mathchardef\nVdash="3\msb@31
\mathchardef\nvDash="3\msb@32
\mathchardef\nVDash="3\msb@33
\mathchardef\ntrianglerighteq="3\msb@34
\mathchardef\ntrianglelefteq="3\msb@35
\mathchardef\ntriangleleft="3\msb@36
\mathchardef\ntriangleright="3\msb@37
\mathchardef\nleftarrow="3\msb@38
\mathchardef\nrightarrow="3\msb@39
\mathchardef\nLeftarrow="3\msb@3A
\mathchardef\nRightarrow="3\msb@3B
\mathchardef\nLeftrightarrow="3\msb@3C
\mathchardef\nleftrightarrow="3\msb@3D
\mathchardef\divideontimes="2\msb@3E
\mathchardef\varnothing="0\msb@3F
\mathchardef\nexists="0\msb@40
\mathchardef\mho="0\msb@66
\mathchardef\thorn="0\msb@67
\mathchardef\beth="0\msb@69
\mathchardef\gimel="0\msb@6A
\mathchardef\daleth="0\msb@6B
\mathchardef\lessdot="3\msb@6C
\mathchardef\gtrdot="3\msb@6D
\mathchardef\ltimes="2\msb@6E
\mathchardef\rtimes="2\msb@6F
\mathchardef\shortmid="3\msb@70
\mathchardef\shortparallel="3\msb@71
\mathchardef\smallsetminus="2\msb@72
\mathchardef\thicksim="3\msb@73
\mathchardef\thickapprox="3\msb@74
\mathchardef\approxeq="3\msb@75
\mathchardef\succapprox="3\msb@76
\mathchardef\precapprox="3\msb@77
\mathchardef\curvearrowleft="3\msb@78
\mathchardef\curvearrowright="3\msb@79
\mathchardef\digamma="0\msb@7A
\mathchardef\varkappa="0\msb@7B
\mathchardef\hslash="0\msb@7D
\mathchardef\hbar="0\msb@7E
\mathchardef\backepsilon="3\msb@7F
\def\Bbb{\ifmmode\let\next\Bbb@\else
 \def\next{\errmessage{Use \string\Bbb\space only in math mode}}\fi\next}
\def\Bbb@#1{{\Bbb@@{#1}}}
\def\Bbb@@#1{\fam\msbfam#1}

\mathchardef\cg="2\eufm@67
\mathchardef\cm="2\eufm@6D

\catcode`\@=\active

\def\eps{{\epsilon}}

\def\<{\langle}
\def\>{\rangle}

\def\tens{\mathop{\otimes}}

\def\ev{{\rm ev}}

\def\id{{\rm id}}

\def\id{{\rm id}}
\def\vec{{\rm Vec}\,}

\def\proof{\goodbreak\noindent{\bf Proof\quad}}

\def\text#1{{\rm #1}}
\def\note#1{}

\documentclass[11pt]{article}
\newtheorem{lemma}{Lemma}[section]
\newtheorem{propos}[lemma]{Proposition}
\newtheorem{example}[lemma]{Example}
\newtheorem{theorem}[lemma]{Theorem}
\newtheorem{cor}[lemma]{Corollary}

\newtheorem{defin}[lemma]{Definition}
\newtheorem{remark}[lemma]{Remark}

\begin{document}
\baselineskip 22pt


\begin{center} {\Large Braiding and exponentiating noncommutative vector
fields}\\
\bigskip
 \large{Edwin J.
Beggs\\{\ }\\Department of Mathematics,\\University of Wales
  Swansea, UK}
\end{center}

{\bf Abstract} The purpose of this paper is to put into a noncommutative
context basic notions related to vector fields from classical
differential geometry.  The manner of exposition is an attempt to make
the material as accessible as possible to classical geometers.  The
definition of vector field used is a specialisation of the Cartan pair
definition, and the paper relies on the idea of generalised braidings
of 1-forms.  The paper considers Kroneker deltas, interior products, Lie
derivatives, Lie brackets, exponentiation of vector fields and
parallel transport.

\vskip0.3in

\section{Introduction}
Classical differential geometry is heavily reliant on 
the use of vector fields, and they also provide
an intuitive way to think about the geometry, linking with ideas of 
flow or motion from physics.  However noncommutative differential
geometry has
been largely concerned with forms.  In this paper I have used
a specialisation of the Cartan pair definition \cite{Borowiec}
of a noncommutative vector field.  My intention was to try to 
formulate noncommutative analogues of certain classical constructions
requiring vector fields, especially the interior product.
Note that the Cartan pair definition of vector fields 
on Hopf algebras was considered in \cite{pjdl}.

The feature which allows us to make any sense of many classical
constructions in the noncommutative world is the generalised `braiding' 
(in some
cases this word is interpreted rather loosely) between bimodules and
1-forms, which is described in \cite{Madore}.  If we consider the braiding in
the commutative case, it is just order reversal of forms or vector
fields.  The fact that this is an honest braiding (i.e.\ satisfies the
braid relation), and that it precisely determines all the differential
forms given just the 1-forms by antisymmetry, becomes the dominant
feature of the commutative case.  Also of great geometrical importance is
the interior product, a pairing between the vector fields and forms
which reduces the degree of the form by one.  A noncommutative
differential calculus with these features (such as the calculus for
the noncommutative torus given in \cite{bdr}) behaves more or less the
same as a commutative differential calculus.

One seemingly strange feature is the number of conditions needed on
the differential calculus for some of the results to hold.  A little
while spent constructing differential calculi on algebras given in
terms of generators and relations will reveal a reason for this. 
There are often a very large number of possible differential calculi
once the constraints of commutativity are removed, and some calculi,
and some covariant derivatives, are nicer than others.  In particular
we arrive at an idea of a compatibility between the differential
calculus and covaraint derivatives and their associated braidings.

The paper begins with standard material \cite{ConnesNCDG}
on differential calculi
 and connections.  Then it considers paired connections and
braidings on vector fields.  The central idea introduced in the paper
is a noncommutative analogue of interior product of a vector field
with an $n$-form.  From here it is not difficult to introduce the Lie
derivative of an $n$-form.  Antisymmetric tensor products of fields
are introduced, and are used to define vector field versions of
curvature and torsion, as well as an idea of Lie bracket.  From a
noncommutative Kroneker delta we can define a differential dimension
of the algebra, which depends on the differential calculus and the
associated braiding.  The example of the noncommutative torus
\cite{bdr} is considered, and proves to be very like the classical
case.  The noncommutative sphere \cite{CHZ} illustrates some rather
less classical behaviour.

The paper ends by considering noncommutative analogues of exponentiation
of vector fields, parallel transport and geodesics.  The problem here is that
the result of an exponentiation is not in general an algebra map. 
However it retains the structure of a cochain map, and is shown to be well
behaved under the coaction of a Hopf algebra on the algebra. 
In the process of doing this we must consider exponentials of `Lie algebra'
elements for the Hopf algebra.  An example of exponentiation is given
on the noncommutative torus.

In the notation, I have made use of overloading certain symbols,
with distinction being made by considering the domains, rather than
have a multiplicity of symbols or indices. 
I use $\id^n$ to be $\id\tens\id\tens\dots\tens\id$ $n$ times. All algebras
are assumed to be unital and associative.

I would like to thank T. Brzezi\'nski (Swansea) and S. Majid (QMW
London) for their help in the preparation of this paper.

\section{Noncommutative differential calculi}

\begin{defin}
A differential structure on an algebra $B$ is a graded algebra
$\Omega^{n}B$ for $n\ge 0$ (i.e.\ there is a multiplication
$\wedge:\Omega^{n}B\tens \Omega^{m}B\to \Omega^{n+m}B$) with
$\Omega^{0}B=B$. In particular the graded algebra structure makes every 
$\Omega^{n}B$ into an
$\Omega^{0}B=B$-bimodule, and we use a dot for this operation, rather
than $\wedge$.  To every differential form $\omega\in\Omega^{n}B$ we
assign a grade $|\omega|=n$.  There is a differential
$d:\Omega^{n}B\to \Omega^{n+1}B$ with $d^{2}=0$ and
$d(\tau\wedge\omega)=d\tau\wedge\omega +(-1)^{|\tau|}\tau\wedge
d\omega$.  In addition we assume that $B.d\Omega^{n}B$ and
$d\Omega^{n}B.B$ are dense in $\Omega^{n+1}B$.
\end{defin}

\begin{defin}
Given differential structures on algebras $B$ and $C$, an algebra map
$f:B\to C$ is called differentiable if there is a well defined map
$f_*:\Omega^1 B\to \Omega^1 C$ defined by $f_*(b\, db')=f(b)\,
df(b')$. Here `well defined' means that if a sum of elements of the form
$b\, db'$ vanishes in $\Omega^1 B$, then the corresponding sum
of $f(b)\,df(b')$ vanishes in $\Omega^1 C$.

  Then $f_*$ is a $B$-$B$ bimodule map, where the left and right
action by $b\in B$ on $\Omega^1 C$ is respectively left and right
multiplication by $f(b)$.  If $g:C\to E$ is also differentiable, then
$g\circ f:B\to E$ is differentiable, and $(g\circ f)_*=g_*\circ
f_*:\Omega^1 B\to\Omega^1 E$.
\end{defin}

\begin{defin} \label{tensdef}
    Given differentiable structures on
algebras $B$ and $C$, the tensor product differential structure on
$B\tens C$ is defined by $\Omega^1(B\tens C)=(\Omega^1B\tens C) \oplus
(B\tens \Omega^1C)$ and $d(b\tens c)=db\tens c+b\tens dc$.  We use
this splitting to define projections $\Pi_1:\Omega^1(B\tens
C)\to\Omega^1B\tens C$ and $\Pi_2:\Omega^1(B\tens C)\to B\tens
\Omega^1C$.  These splittings obey the functorial conditions that if
$f:B\to M$ and $g:C\to A$ are differentiable algebra maps, then
$\Pi_1\circ(f\tens g)_{*}=(f_{*}\tens g)\circ\Pi_1:\Omega^1(B\tens C)\to
\Omega^{1}M\tens
A$ and $\Pi_2\circ(f\tens g)_{*}=(f\tens g_{*})\circ\Pi_2:\Omega^1(B\tens C)\to
M\tens\Omega^{1}A$.
\end{defin}

\section{Covariant derivatives and braidings with 1-forms}
We shall take $M$ to be an algebra with a specified differential calculus. 
If $M$ were the algebra of functions on a topological space,
then given a bundle over the space, the sections of the bundle form 
a module. 
 In the noncommutative setting, we consider modules in place of bundles.

\begin{defin} Given a left $M$-module $E$, a left $M$-covariant derivative
is a map $\nabla:E\to\Omega^{1}M\tens_{M}E$ which obeys the condition
$\nabla(m.e)=dm\tens e+m.\nabla e$
for all $e\in E$ and $m\in M$.
\end{defin}

\begin{remark} The tensor product over $M$, $F\tens_M E$ for a left $M$-module $E$ and a right
$M$-module $F$, is like the usual tensor product of vector spaces, but in addition
we make the identification $f.m\tens e=f\tens m.e$ for all $f\in F$, $e\in E$ and $m\in M$.
\end{remark}

\begin{defin} \label{ppll} (See \cite{Madore}.) A bimodule covariant derivative on an
 $M$-bimodule $E$ is a pair $(\nabla,\sigma)$,
where $\nabla:E\to\Omega^{1}M\tens_{M}E$ 
is a left $M$-covariant derivative, and $\sigma:E\tens_M\Omega^1 M\to 
\Omega^1 M \tens_M E$ is a bimodule map
called the `braiding' (even when it isn't) obeying
\[
\nabla(e.m)\,=\,\nabla(e).m\,+\,\sigma(e\tens dm)\ .
\]
\end{defin}

\begin{remark}\label{rem1345} Of course, given a left $M$-covariant derivative
$\nabla$ on a bimodule $E$, we can try to define a compatible left braiding
by $\sigma(e\tens a\,db)=\nabla(e.ab)-\nabla(e.a).b$.  However this
might not give a well defined result,
but we do see that there is at most one braiding compatible with a
given connection.  From this formula we also see (using the fact that
$d$ is a derivation) that the resulting $\sigma$ (if well defined) is an $M$-bimodule map.
This means, as we are only concerned with braidings compatible with
connections, that there is no point in weakening the definition of a
braiding to a left or right module map.
\end{remark}

\begin{propos} Given $(\nabla',\sigma')$ a bimodule covariant derivative
 on the bimodule $E$, any other left covariant derivative
 on the bimodule $E$ is of the form $\nabla=\nabla'+\Gamma$,
where $\Gamma:E\to
\Omega^{1}M\tens_{M}E$ is a left $M$-module map. We get $(\nabla,\sigma)$ a bimodule covariant derivative
if and only if the braiding
 $\sigma(e\tens a\,db)=\sigma'(e\tens a\,db)+
\Gamma(e.ab)-\Gamma(e.a).b$ is well defined. 
In particular $\sigma=\sigma'$ if and only if
$\Gamma$ is an $M$-bimodule map.
\end{propos}
\proof Straightforward.
\quad$\square$

\begin{propos} \label{tens45} (See \cite{Madore}.)  Given $(\nabla,\sigma_E)$ a
bimodule covariant derivative on the bimodule $E$ and $\nabla$ a left
covariant derivative on the left module $F$, there is a left
$M$-covariant derivative on $E\tens_{M}F$ given by $\nabla\tens{\rm
id}_{F}+(\sigma_E\tens{\rm id}_{F})({\rm id}_{E}\tens\nabla)$. 
Further if $F$ is also an $M$-bimodule with a bimodule covariant
derivative $(\nabla,\sigma_F)$, then there is a compatible braiding on
$E\tens_{M}F$ given by $ \sigma_{E\tens F} \,=\, (\sigma_{E}\tens{\rm
id})(\id\tens \sigma_{F}) $.
\end{propos}
\proof Applying the formula to $e\tens m.f$ we get $\nabla e\tens
m.f+(\sigma_E\tens{\rm id}_{F})(e\tens dm\tens f+e\tens m.\nabla f)$. 
Aplying the formula to $e.m\tens f$ we get $\nabla(e.m)\tens f+
(\sigma_E\tens{\rm id}_{F})(e.m\tens \nabla f)$, and these are the same by
definition of $\sigma_E$.
 This shows that the given formula is
well defined on $E\tens_{M}F$.  The left multiplication property is
true because $\sigma_E$ is a left $M$-module map.

For the second part, we use the formula from \ref{rem1345},
\begin{eqnarray*}
    \sigma_{E\tens F}(e\tens f\tens a.db) &=& \nabla(e\tens f.ab)-
    \nabla(e\tens f.a).b \cr
    &=& \nabla(e)\tens f.ab + (\sigma_{E}\tens{\rm id})(e\tens \nabla(f.ab)) \cr
    &&- \nabla(e)\tens f.ab-(\sigma_{E}\tens{\rm id})(e\tens
    \nabla(f.a)b)\ .\quad\square
\end{eqnarray*}

\begin{defin}\label{pres3}
A left module map $\theta:E\to F$ is said to be preserved by the covariant derivatives
$\nabla$ on $E$ and $F$ if $\nabla\circ\theta=(\id\tens\theta)\nabla:E\to\Omega^1 M\tens_M F$.
\end{defin}

\begin{propos} Given bimodules $E$ and $F$ with left covariant derivatives 
$(\nabla,\sigma)$,
the bimodule map $\theta:E\to F$ obeys the condition $(\id \tens\theta)
    \sigma=\sigma(\theta\tens\id)$ if and only if the map
    $\nabla\circ\theta-(\id\tens\theta)\nabla :E\to \Omega^{1}M\tens_{M}
    F$ is an $M$-bimodule map.
\end{propos}
\proof The left module map property of $\nabla\theta-(\id\tens\theta)\nabla$ is fairly simple.
The right module map property is given by subtracting the following equations:
\begin{eqnarray*}
(\id \tens\theta)\sigma(e\tens db) &=& (\id
\tens\theta)(\nabla(e.b)-\nabla(e).b)\ ,\cr \sigma(\theta (e)\tens db) &=&
\nabla(\theta (e).b)-\nabla(\theta (e)).b\ .\quad\square
\end{eqnarray*}

\begin{example} The simplest $M$-bimodule is $M$ itself. Unless otherwise stated, 
we take the covariant derivative $\nabla=d:M\to\Omega^1 M\tens_M M=\Omega^1 M$. The
corresponding $\sigma$ is the identity.
\end{example}

\begin{defin}\label{derres} Given a left covariant derivative $\nabla$ on a
left $M$-module $F$ and a left submodule $G\subset F$, $\nabla$ is
said to restrict to $G$ if $\nabla G\subset \Omega^{1}M\tens_M G$.
If a covariant derivative preserves a
sub-bimodule, then its associated braiding also preserves the sub-bimodule,
i.e.\ $\sigma(G\tens_M\Omega^1 M)\subset \Omega^1 M\tens_M G$.
\end{defin}

\section{Finitely generated projective modules}
General modules over algebras can be quite badly behaved, so here
we offer a definition and some results about a
well known nice class of modules, the
finitely generated projective modules. See \cite{fgp} for more details.

\begin{defin} The dual $E^*$ of a right 
$M$-module $E$ is defined to be $Hom_M(E,M)$, the right module maps from $E$ to $M$.
Then $E^*$ has a left module structure given by $(m.\alpha)(e)=m.\alpha(e)$ for all
 $\alpha\in E^*$ and $e\in E$. 
If $E$ is a bimodule, then $E^*$ has 
 a right module structure given by $(\alpha.m)(e)=\alpha(m.e)$, and there is a bimodule map evaluation
$\ev:E^*\tens_M E\to M$. 
\end{defin}

\begin{defin}\label{canon} A right $M$-module $E$ is said to be finitely 
    generated projective if
there are $e_i\in E$ and $\alpha_i\in E^*$ (for integer $1\le i\le n$)
(the `dual basis') so that for all $e\in E$, $e=\sum e_i.\alpha_i(e)$. 
From this it follows directly that $\alpha=\sum \alpha(e_i).\alpha_i$
for all $\alpha\in E^*$.
\end{defin}

\begin{example}\label{classdual}
This condition may seem rather esoteric, but it has a simple example. In $C^\infty(\Bbb R^n)$ with 
coordinates 
$\{x^1,\dots,x^n\}$ the sections $\Omega^1(\Bbb R^n)$ of the cotangent $T^*\Bbb R^n$ bundle 
has a module basis $dx^1\dots dx^n$ (i.e.\ 
every section of $T^*\Bbb R^n$ can be written as a sum of functions times the basis elements). 
The dual basis is $\partial_j\in \Omega^1(\Bbb R^n)^*$
($1\le j\le n$) where $\partial_j(dx^i)=\delta^i_j$. The reader should note that the dual basis
$(dx^i,\partial_i)$ is definitely not unique, though we will see shortly that a unique object can be made
by combining them. 
Classically the 
complication comes
when considering a manifold made by patching together coordinate charts. Then we have to apply partitions of unity to
the previous construction on each coordinate chart. Of course, the dual of the 1-forms is the vector fields, 
but we should save that fact for later. 
\end{example}

\begin{propos}\label{endo} If an $M$-bimodule $E$ is finitely generated projective,
and $F$ is a right $M$-module,
 there is an isomorphism $\vartheta:F\tens_M E^* \to {\rm Hom}_M(E,F)$
defined by $\vartheta(f\tens\alpha)(e)=f.\alpha(e)$.  
\end{propos}
\proof The inverse map is $\vartheta^{-1}(T)=\sum T(e_i)\tens\alpha_i$.\quad$\square$

\begin{propos}\label{endol} If an $M$-bimodule $E$ is finitely generated projective,
and $F$ is a left $M$-module,
 there is an isomorphism $\varphi:E\tens_M F \to {}_M{\rm Hom}(E^*,F)$
(the left module maps from $E^*$ to $F$)
defined by $\varphi(e\tens f)(\alpha)=\alpha(e).f$.  
\end{propos}
\proof The inverse map is $\varphi^{-1}(T)=\sum e_i\tens T(\alpha_i)$.\quad$\square$

\begin{cor}\label{corsub} Suppose that we have a map $T:F\to H$ between left
$M$-modules, with kernel $K\subset F$. Then
for a finitely generated projective $M$-bimodule $E$, the map $\id\tens T:E\tens_M F\to E\tens_M H$
has kernel $E\tens_M K$.   
\end{cor}
\proof If we use the isomorphism in \ref{endol}, we get the map
$T\circ:{}_M{\rm Hom}(E^*,F)\to {}_M{\rm Hom}(E^*,H)$, and this has kernel
${}_M{\rm Hom}(E^*,K)$.\quad$\square$

\section{Evaluations and coevaluations}
From now on, we take all right modules considered in the paper to be 
{\it finitely generated projective}.
Also suppose that the bimodules have a bimodule covariant derivative $(\nabla,\sigma)$,
and that $\sigma$ is invertible. 

\begin{propos} \label{camab} Given a bimodule covariant derivative
$(\nabla,\sigma_E)$ on the $M$-bimodule $E$ for which the braiding is
invertible, there is a unique bimodule covariant derivative
$(\nabla,\sigma_{E^*})$ on $E^*$ so that the map $\ev:E^*\tens_M E\to
M$ is preserved by the covariant derivatives (see \ref{pres3}).  It is
defined in terms of the dual basis $(e_i,\alpha_i)$ of $E$ given in
\ref{canon} by
\begin{eqnarray*}
\sigma_{E^*}(\alpha\tens\xi) &=& \sum (\ev\tens\id)(\id\tens\sigma_E^{-1})(\alpha\tens\xi\tens e_i) \tens \alpha_i\ ,
\cr
\nabla\alpha &=& \sum d(\alpha(e_i))\tens\alpha_i\,-\,\sum(\ev\tens\id)(\id\tens\sigma_E^{-1})(\id\tens\nabla)(\alpha\tens e_i)
\tens\alpha_i\ .
\end{eqnarray*}
\end{propos}
\proof First we check that the formulae give a left covariant derivative:
\begin{eqnarray*}
\nabla(m.\alpha) &=& \sum d(m.\alpha(e_i))\tens\alpha_i\,-\,\sum(\ev\tens\id)
(\id\tens\sigma_E^{-1})(\id\tens\nabla)(m.\alpha\tens e_i)
\tens\alpha_i \cr
&=& dm\tens\sum \alpha(e_i).\alpha_i \,+\, m.\nabla(\alpha)\,=\,
dm\tens \alpha \,+\, m.\nabla(\alpha)\ .
\end{eqnarray*}
 Given that the braiding on $M$ is trivial, the condition that $\ev:E^*\tens_M E\to M$ preserves the braiding is
\begin{eqnarray}\label{k}
\ev\tens\id=(\id\tens\ev)\sigma_{E^*\tens E}=(\id\tens\ev)(\sigma_{E^*}\tens\id)(\id\tens\sigma_E):
E^*\tens_M E\tens_M\Omega^1M\to \Omega^1 M\ .
\end{eqnarray}
We check this by
\begin{eqnarray}\label{brad-eval}
(\id\tens\ev)(\sigma_{E^*}\tens \id)(\alpha\tens\xi\tens e) &=& 
\sum (\id\tens\ev)\Big(
(\ev\tens\id)(\id\tens\sigma_E^{-1})(\alpha\tens\xi\tens e_i) \tens \alpha_i\tens e\Big) \cr
&=& \sum
(\ev\tens\id)(\id\tens\sigma_E^{-1})(\alpha\tens\xi\tens e_i.\alpha_i(e)) \cr
&=& 
(\ev\tens\id)(\id\tens\sigma_E^{-1})(\alpha\tens\xi\tens e)\ .
\end{eqnarray}
The $\sigma_{E^*}$ with this property (\ref{brad-eval}) is unique by \ref{endo}.
To see that $\nabla$ preserves the evaluation:
\begin{eqnarray}\label{nabla-eval}
(\id\tens\ev)(\nabla(\alpha)\tens e) &=& 
\sum d(\alpha(e_i)).\alpha_i(e)\,-\,\sum(\ev\tens\id)(\id\tens\sigma_E^{-1})(\id\tens\nabla)
(\alpha\tens e_i)
.\alpha_i(e) \cr
&=&  d(\alpha(e))
\,-\,(\ev\tens\id)(\id\tens\sigma_E^{-1})(\id\tens\nabla)
(\alpha\tens e) \cr
&& -\sum \alpha(e_i)\,d(\alpha_i(e))\,+\,\sum(\ev\tens\id)(\alpha\tens e_i\tens d(\alpha_i(e)))\cr 
&=&
d(\alpha(e))\,-\, (\id\tens\ev)(\sigma_{E^*}\tens\id)(\id\tens\nabla)(\alpha\tens e)\ .
\end{eqnarray}
The $\nabla$ with this property (\ref{nabla-eval}) is unique by \ref{endo}.
Finally we check the compatibility condition in \ref{ppll}, using (\ref{nabla-eval}):
\begin{eqnarray*}
(\id\tens\ev)(\nabla(\alpha.m)\tens e) &=&
  d((\alpha.m)(e))
\,-\,(\ev\tens\id)(\id\tens\sigma_E^{-1})(\id\tens\nabla)(\alpha.m\tens e) \cr
 &=&
  d(\alpha(m.e))
\,-\,(\ev\tens\id)(\id\tens\sigma_E^{-1})(\id\tens\nabla)(\alpha\tens m.e) \cr
&&+\,(\ev\tens\id)(\id\tens\sigma_E^{-1})(\alpha\tens dm\tens e) \cr
&=& (\id\tens\ev)(\nabla(\alpha).m\tens e)\,+\, (\id\tens\ev)(\sigma_{E^*}\tens\id)(\alpha\tens dm\tens e)\ .
\quad\square
\end{eqnarray*}

\begin{defin}\label{kroneker} Given an $M$-bimodule $E$, the Kroneker delta $\delta_E\in E\tens_M E^*$ is defined so that
$(\id_E\tens \ev)(\delta_E\tens e)=e$ for all $e\in E$. In terms of tensor categories,
$\delta_E$ is a coevaluation. For $E$ finitely generated projective
(see \ref{canon}), we have $\delta_E=\sum e_i\tens\alpha_i$.
\end{defin}

\begin{propos}\label{kronprop} For a $\delta_E$ given in \ref{kroneker}:

a)\quad $\delta_E$ is unique. 

b)\quad $m.\delta_E=\delta_E.m$ for all $m\in M$. 

c)\quad $\sigma_{E\tens_M E^*}(\delta_E\tens\xi)=\xi\tens\delta_E$ for all $\xi\in\Omega^1 M$.

d)\quad $\nabla(\delta_E)=0$. 
\end{propos}
\proof $\delta_E$ corresponds to the identity map under the isomorphism in \ref{endo},
proving (a). By \ref{endo} again, to prove (b) we only have to show that
\begin{eqnarray}\label{test1}
m.e &=& (\id\tens\ev)(\delta_E.m\tens e)\ .
\end{eqnarray}
But the right hand side of (\ref{test1}) is
\[
\sum e_i.(\alpha_i.m)(e)\,=\,\sum e_i.(\alpha_i)(m.e)\,=\,m.e
\]
for all $e\in E$ as required. 
By \ref{endo} again 
and using the fact that $\sigma_E$ is invertible,
 to prove (c) we only have to show that, for all $\xi\in\Omega^1 M$ and $e\in E$,
\begin{eqnarray}\label{test2}
\sigma_E(e\tens\xi) &=& (\id^2\tens\ev)(\sigma_{E\tens E^*}\tens\id)(\delta_E\tens\sigma_E(e\tens\xi))\ .
\end{eqnarray}
From (\ref{k}), the right hand side of (\ref{test2}) is
\begin{eqnarray*}
(\id^2\tens\ev)(\sigma_{E}\tens\id^2)
(\id\tens\sigma_{E^*}\tens\id)(\delta_E\tens\sigma_E(e\tens\xi)) &=& \sigma_E(\id\tens\ev\tens\id)
(\delta_E\tens e\tens\xi)\cr
&=& \sigma_E(e\tens \xi)\ ,
\end{eqnarray*}
as required. By \ref{endo} again, to prove (d) we only have to show that
$(\id^2\tens\ev)(\nabla\delta_E\tens e)=0$ for all $e\in E$. Then, using (\ref{nabla-eval}),
\begin{eqnarray*}
(\id^2\tens\ev)(\nabla\delta_E\tens e) &=& (\id^2\tens\ev)\sum\Big(\nabla e_i\tens\alpha_i\tens e
+(\sigma_E\tens\id^2)(e_i\tens\nabla\alpha_i\tens e)\Big) \cr
&=& \sum\Big((\nabla e_i).\alpha_i(e)\,+\,\sigma_E(e_i\tens(\id\tens\ev)(\nabla\alpha_i\tens e))\Big)\cr
&=&\sum\Big( (\nabla e_i).\alpha_i(e) \,+\, \sigma_E(e_i\tens d\alpha_i(e))\cr &&-\,
\sigma_E(e_i\tens(\id\tens\ev)(\sigma_{E^*}\tens\id)(\alpha_i\tens\nabla e))\Big) \cr
&=&\sum\Big( \nabla(e_i.\alpha_i(e))\,-\, \sigma_E(e_i\tens(\id\tens\ev)(\sigma_{E^*}\tens\id)(\alpha_i\tens\nabla e))\Big) \cr
&=& \nabla(e)\,-\,\sum \sigma_E(e_i\tens(\ev\tens\id)(\id\tens\sigma_{E}^{-1})(\alpha_i\tens\nabla e))\ .
\end{eqnarray*}
Now substitute $\sigma_{E}^{-1}\nabla e=\sum f_j\tens\eta_j\in E\tens\Omega^1 M$, giving
\begin{eqnarray*}
(\id^2\tens\ev)(\nabla\delta_E\tens e)&=& \nabla(e)\,-\,\sum 
\sigma_E(e_i\tens \alpha_i(f_j).\eta_j  ) \cr
&=& \nabla(e)\,-\,\sum 
\sigma_E(e_i.\alpha_i(f_j)\tens\eta_j  ) \cr
&=& \nabla(e)\,-\,\sum 
\sigma_E(f_j\tens\eta_j  ) \,=\,0\ .\quad\square
\end{eqnarray*}

\section{Vector fields}
In this section we assume that the $M$-bimodule $\Omega^1 M$ is finitely generated projective
as a right module. 
We use the dual basis $\xi_i\in \Omega^1 M$ and $X_i\in (\Omega^1M)^*$ so that $\sum
\xi_i.X_i(\eta)=\eta$ for all $\eta\in \Omega^1 M$.

\begin{defin} 
Define the vector fields on an algebra $M$ by ${\rm Vec}\,M=(\Omega^1M)^*$,
the right $M$-module maps from $\Omega^{1}M$ to $M$. Then evaluation gives a
bimodule map $\ev:\vec M\tens_{M}\Omega^{1}M\to M$. If $f:B\to M$ is a differentiable
 algebra map, we define
   $f^{*}:{\rm Vec}\,M\to {\rm Vec}\,B$ by $(f^{*}X)(\xi)=X(f_{*}\xi)$
   for $\xi\in\Omega^{1}B$.
\end{defin}

\begin{defin}
 An $X\in {\rm Vec}\,M$ gives a `directional derivative' map
 $D_X:M\to M$ defined by
\[
{\rm Vec}\,M \tens M \stackrel{ {\rm id}\tens d}{\longrightarrow}
{\rm Vec}\,M \tens \Omega^1 M
\stackrel{ {\rm eval}}{\longrightarrow} M\ .
\]
This map is a derivation on $M$ if and only if $X:\Omega^1 M\to M$ 
is also a left $M$-module map. 
In general $D_{X}(a\,b)=D_{X}(a)\,b\,+\,D_{X. a}(b)$.

If the left $M$-module $E$ has a left $M$-covariant derivative $\nabla$, then
given $X\in\vec M$ we define the covariant directional derivative by
$\nabla_{X}e=(\ev\tens \id)(X\tens\nabla e)$.  The reason for defining
vector fields as right $M$-module maps was so that this would be well defined.
\end{defin}

\begin{defin}\label{vecbraid} The vector fields are braided by 
    $\sigma^{-1}:\vec M\tens \vec M\to \vec M
\tens\vec M$ given by the formula
$\sigma^{-1}(X\tens Y)=\sum (\ev\tens\id)(X\tens\sigma_{\vec M}
(Y\tens\xi_i))\tens X_i$. Note that
we use the notation $\sigma^{-1}$ to fit with the crossings in a braided
category - we do not claim that $(\sigma^{-1})^{-1}$ exists!
\end{defin}

\begin{propos}\label{dubrad} The braiding in \ref{vecbraid} is the unique
braiding for which
\[
(\id\tens\ev)(\sigma^{-1}\tens\id) \,=\, (\ev\tens\id)(\id\tens\sigma_{\vec M}):
\vec M\tens_M\vec M\tens_M \Omega^1 M\to \vec M\ .
\]
\end{propos}
\proof
Uniqueness follows from \ref{endo} again. 
For all $\eta\in \Omega^1
M$ and $X,Y\in \vec M$,
\begin{eqnarray*}
(\id\tens\ev)(\sigma^{-1}\tens\id)(X\tens Y\tens \eta) &=& \sum 
(\ev\tens\id)(X\tens\sigma_{\vec M}(Y\tens\xi_i)). X_i(\eta) \cr
&=& (\ev\tens\id)(X\tens\sigma_{\vec M}(Y\tens\eta))
\ .\quad \square
\end{eqnarray*}

\begin{defin} Using the fact that $\Omega^1 M$ is finitely generated projective,
we have a unique Kroneker delta (see \ref{kroneker})
$\delta\in \Omega^1 M\tens_M\vec M$. In addition we define 
$\hat\delta=\sigma^{-1}\delta\in\vec M\tens_M \Omega^1 M$,
and $\dim M=\ev(\hat\delta)\in M$. Note that $\dim M$ is a 
central element in $M$ by \ref{kronprop}(b).
\end{defin}

\begin{remark} Following from \ref{classdual},
note that $\delta$ and $\hat\delta$ in classical differential geometry are just 
the usual Kroneker deltas, $\delta^i_{\phantom{i}j}$ and $\delta_i^{\phantom{i}j}$.
It is then immediate that $\dim M$ is a constant function with value 
the dimension of the manifold. 
\end{remark}

\section{Interior products} \label{sectinprod}
In this section we would like to define the interior product of a
vector field with an $n$-form.  However we must remember that the 
$n$-forms are not realised as a subspace of the $n$-fold tensor
product of the 1-forms, but rather as a quotient of them
by $\Theta^n M=\ker\wedge:\bigotimes_M^n \Omega^1 M\to \Omega^n M$.  This leads
us to a compatibility condition between the braiding and the
differential calculus which is necessary to define interior products
with $\Omega^{n}M$.

\begin{defin} Recursively define the $M$-bimodule map $
\sigma_n:\bigotimes_M^n \Omega^1 M\to \bigotimes_M^n \Omega^1 M$,
beginning with $\sigma_1=\id$ and $\sigma_2=\sigma$, and continuing
with $\sigma_{n+1} \,=\, (\sigma\tens\id^{n-1})(\id\tens\sigma_n)$. 
It is easy to see that for all $r,s\ge 0$,
$(\sigma_{s+1}\tens\id^r)(\id^s\tens\sigma_{r+1})=\sigma_{r+s+1}$.
\end{defin}

\begin{defin} \label{intdef} Define the interior product $X\lrcorner\,
z\in\bigotimes_M^{n-1}\Omega^{1}M$ for $X\in\vec M$ and
$z\in\bigotimes_M^{n}\Omega^{1}M$ as $(\ev\tens\id^{n-1})(X\tens
T_n(z))$, where
\[
T_n\,=\,-\,\sum_{r=1}^n (-1)^r\,\sigma_r\tens \id^{n-r} :
\bigotimes_M^{n}\Omega^{1}M \to \bigotimes_M^{n}\Omega^{1}M\ .
\]
For $\omega\in\Omega^1 M$ we have $X\lrcorner\,\omega=X(\omega)$.
\end{defin}

\begin{propos}
 The map $\lrcorner:\vec
    M\tens_{M}(\tens^{n+1}_{M}\Omega^{1}M)\to \tens^{n}_{M}\Omega^{1}M$ is an $M$-bimodule map. 
\end{propos}
\proof All its component maps are bimodule maps.\quad$\square$

\begin{defin}\label{compatform} The interior product operation is said to
be compatible with the differential calculus if
 $T_{n+1}(
\Theta^{n+1}M)\subset \Omega^1 M\tens_M\Theta^{n}M$ for all $n\ge 1$.
In this case, we get an interior product $\lrcorner:\vec
M\tens_{M}\Omega^{n+1}M\to\Omega^{n}M$. We conventionally add $X\lrcorner\,m=0$ for $m\in\Omega^0 M$.
\end{defin}

\begin{propos}\label{lrcomp} If $\lrcorner$ is compatible with the 
    differential calculus, then

a)\quad $\sigma_{n+1}(\Theta^n M\tens_M \Omega^1
M)\subset \Omega^1 M\tens_M \Theta^n M$ for all $n\ge 1$. 

b)\quad $\Theta^{2}M$ is contained
in the +1 eigenspace of $\sigma:
\Omega^{1}M\tens_{M}\Omega^{1}M \to \Omega^{1}M\tens_{M}\Omega^{1}M $.
\end{propos}
\proof To prove (a), given $z\in \Theta^n M$ and $\xi\in\Omega^1 M$,
we know that $z\tens\xi\in \Theta^{n+1} M$, so
$T_{n+1}(z\tens\xi)\in \Omega^1 M\tens_M \Theta^n M$ by our assumption. But
\[
T_{n+1}(z\tens\xi)\,=\,T_n(z)\tens\xi\,+\,(-1)^n \sigma_{n+1}(z\tens\xi)\ ,
\]
and, also by our assumption, $T_n(z)\tens\xi\in \Omega^1 M\tens \Theta^{n-1}M\tens\Omega^1 M\subset
\Omega^1 M\tens_M \Theta^n M$. We deduce that $\sigma_{n+1}(z\tens\xi)\in \Omega^1 M\tens_M \Theta^n M$.

To prove (b), note that $\Theta^1M=0$, so we have $T_2 \Theta^2 M=0$ by our assumption.
 \quad$\square$

\section{Lie derivatives of forms}
Having defined interior products of vector fields with $\Omega^{n}M$
in section \ref{sectinprod}, we are in the happy position of being able to
define the Lie derivative of an $n$-form with respect to a vector
field.  We assume that we have
covariant derivatives and braidings satisfying \ref{compatform}.

\begin{defin} \label{liederdef} We define the Lie derivative of
$\omega\in\Omega^{n}M$ with respect to $X\in{\rm Vec}\,M$ to be ${\cal
L}_{X}\omega=X\lrcorner(d\omega)+ d(X\lrcorner\, \omega)$.  Note that
${\cal L}_{X}(m)=D_{X}m$ for $m\in\Omega^{0}M=M$.
\end{defin}

\begin{propos} \label{lieprops} The Lie derivative ${\cal L}:\vec M\tens
\Omega^{n}M\to \Omega^{n}M$ obeys the following rules:

a)\quad
${\cal L}_X(m.\omega)={\cal L}_{X.m}(\omega)+X\lrcorner(dm\wedge \omega)$.

b)\quad
${\cal L}_{m.X}(\omega) \,=\, dm \wedge (X\lrcorner \,\omega)\,+\,m.{\cal
    L}_{X}\omega$.

c)\quad ${\cal L}_X(\omega.m)={\cal L}_X(\omega).m+(-1)^{|\omega|}\,
(X\lrcorner(\omega\wedge dm)-(X\lrcorner\omega)\wedge dm)$.

d)\quad $d\,{\cal L}_{X}={\cal L}_{X}\,d:\Omega^{n}M\to
    \Omega^{n+1}M$.
\end{propos}
\proof More or less immediate from the definition. \quad$\square$

\section{Covariant derivatives of higher forms}

\begin{remark} Given a bimodule covariant derivative $(\nabla,\sigma)$ on
$\Omega^1 M$, the discussion in \ref{tens45} gives a covariant derivative
$(\nabla,\sigma_{n+1})$ on $\bigotimes_M^n \Omega^1 M$, given by
\[
\nabla\,=\,\sum_{i=1}^{n} (\sigma_i\tens \id^{n+1-i})(\id^{i-1}\tens\nabla\tens\id^{n-i})
:\bigotimes_M^n \Omega^1 M\to \Omega^1 M\tens_M\Big(\bigotimes_M^n \Omega^1 M\Big)\ .
\]
\end{remark}

\begin{propos} If $\nabla$ on $\bigotimes_M^n \Omega^1 M$ preserves (in
the sense of \ref{derres}) the submodule $\Theta^n M=\ker\wedge
:\bigotimes_M^n \Omega^1 M\to\Omega^n M$, then we get a covariant
derivative on $\Omega^n M$ by quotienting.
\end{propos}
\proof Reasonably direct from the previous statements.  \quad$\square$

\section{Antisymmetry and Lie brackets of vector fields}

\begin{defin} \label{antidef} An $x\in\vec M\tens\vec M$ is called
antisymmetric if $\ev(\id\tens\ev\tens\id)(\pi x\tens k)=0$ for all
$k\in \Theta^{2}M$, where $\pi$ is the quotient map from $\vec
M\tens\vec M$ to $\vec M\tens_{M}\vec M$.  We call $A^{2}M$ the set of
antisymmetric elements in $\vec M\tens\vec M$.
\end{defin}

\begin{remark} The map $(\ev\tens \ev)(\id\tens\sigma\tens\id): \vec
M\tens_{M} \vec M\tens_{M} \Omega^{1}M \tens_{M} \Omega^{1}M\to M$ can
also be written as
$\ev(\id\tens\ev\tens\id)(\sigma^{-1}\tens\id\tens\id)$ and as
$\ev(\id\tens\ev\tens\id)(\id\tens\id\tens\sigma^{-1})$.  If \ref{compatform} holds,
 it follows that all eigenspaces of
$\sigma^{-1}$ except the +1 eigenspace are contained in $\pi A^2 M$.
\end{remark}

\begin{defin}\label{deflie} Define a map $\phi:A^2M\to {\rm Vec}\,M$
by the following formula, where $\xi\in\Omega^{1}M$:
\[
\phi(X\tens Y)(\xi)\,=\,D_X(Y(\xi))\,+\,\ev(\id\tens \ev\tens\id) (X\tens
Y\tens z)\ ,
\]
where $z\in \Omega^1 M
\tens_M\Omega^1 M$ is chosen so that $\wedge z=d\xi$.  By the previous
discussion the choice does not matter.
 To check that its image is in $\vec M$ we use the
following proposition.
\end{defin}

\begin{propos}\label{prop17} The image of the map $\phi$ in \ref{deflie} is in
    $\vec M$. Further $\phi$ is a left
$M$-module map, but not in general a right module map,
as $\phi(X\tens Y).m=\phi(X\tens Y.m)+X.D_Y(m)$. 
Also $\phi(X\tens m.Y)=\phi(X.m\tens Y)+D_{X}(m).Y$.
\end{propos}
\proof To see that $\phi(X\tens Y)$ is a right module map use the following,
where $\wedge z=d\xi$,
\begin{eqnarray*} 
    \phi(X\tens Y)(\xi.m) &=& D_X(Y(\xi).m)\,+\,\ev(\id\tens
\ev\tens\id) (X\tens Y\tens (z.m-\xi\tens dm)) \cr
&=& \phi(X\tens Y)(\xi).m\,+\, X(Y(\xi).dm)\,-\, X(Y(\xi).dm)\ .
\end{eqnarray*}
It is quite easy to see that $\phi(m.X\tens Y)(\xi)=m.\phi(X\tens Y)(\xi)$. 
For the right action,
\begin{eqnarray*} 
\phi(X\tens Y)(m.\xi) &=& D_X(Y(m.\xi))\,+\,\ev(\id\tens \ev\tens\id) (X\tens
Y\tens (m.z+dm\tens\xi)) \cr
&=& \phi(X\tens Y.m)(\xi)\,+\,\ev(\id\tens \ev\tens\id) (X\tens
Y\tens dm\tens\xi)\ .
\end{eqnarray*}
Finally we calculate
\begin{eqnarray*} 
\phi(X\tens m.Y)(\xi) &=& 
D_X(m.Y(\xi))\,+\,\ev(\id\tens \ev\tens\id) (X\tens
m.Y\tens z) \cr &=& D_X(m).Y(\xi)+
D_{X.m}(Y(\xi))\,+\,\ev(\id\tens \ev\tens\id) (X.m\tens Y\tens z)\
.\quad\square
\end{eqnarray*}

\begin{remark}
 Now $\phi$ is the non-commutative version of the Lie bracket, but it needs
    to be applied to elements of $\vec M\tens \vec M$ which are already
    antisymmetric.  In the commutative case we would have $\phi(X\tens
    Y-Y\tens X)=[X,Y]$.  In the non-commutative case the problem is the
    anologue of the $X\tens Y\mapsto X\tens
    Y-Y\tens X$ operation, i.e.\ how to antisymmetrise elements.

    If $\sigma^{2}=\id$, then we can antisymmetrise on $\vec M\tens_{M}
    \vec M$ by $X\tens Y\mapsto X\tens Y-\sigma(X\tens Y)$, but
$\phi$ does not descend to a well defined map on a
subspace of ${\rm Vec}\,M\tens_M {\rm Vec}\,M$.
 Any definition of Lie bracket will depend an being given an
    antisymmetrisation operation on $\vec M\tens \vec M$. 
\end{remark}

\section{Commutative differential geometry}
In this case $\sigma:\Omega^1 M 
\tens_M \Omega^1 M  \to \Omega^1 M  \tens_M \Omega^1 M $ is just transposition
$\sigma(\xi\tens\eta)=\eta\tens\xi$,
and satisfies the braid relation.  We have $\sigma^2=1$, and $\Omega^1 M
\tens_M \Omega^1 M $ splits into a direct sum of the +1 and -1
eigenspaces of $\sigma$.  The +1 eigenspace is $\Theta^2 M=\ker\wedge:
\Omega^1 M \tens_M \Omega^1 M \to\Omega^2 M$, and we can use this to
identify $\Omega^2 M$ with the -1 eigenspace (antisymmetric tensors)
in $\Omega^1 M \tens_M \Omega^1 M $.

\section{The non-commutative torus}
    Take the algebra $\Bbb T^2_q$ generated by invertible elements $u$ and $v$, subject to 
$uv=qvu$, where $q$ is a unit norm complex number. 
The simplest differential 
calculus (there are many to choose from)
 on $\Bbb T^2_q$ \cite{bdr} is generated by $\{u,v,du,dv\}$, subject to the relations
\begin{eqnarray} \label{nctrel}
du\wedge dv=-q\, dv\wedge du\ &,& u\, dv=q\, dv\, u\ ,\quad v\,
du=q^{-1}du\, v\ ,\cr [u,du]=[v,dv]=0\ &,& du\wedge du=dv\wedge dv=0\ .
\end{eqnarray} 
Then ${\rm Vec}\,\Bbb T^2_q$ is generated as a left $\Bbb T^2_q$ module by
the elements $\partial_{u}$ and $\partial_{v}$, where
$\partial_{u}\,(du)=\partial_{v}\,(dv)=1$ and
$\partial_{u}\,(dv)=\partial_{v}\,(du)=0$.  In fact for any $X\in {\rm
Vec}\,\Bbb T^2_q$, $X=X(du)\,\partial_{u}+X(dv)\,\partial_{v}$.  Now we find
the right actions
\begin{eqnarray}
(X. u)(du)=X(du)\,u\ &,& (X. u)(dv)=q\,X(dv)\,u\ , \cr (X.
v)(du)=q^{-1}\,X(du)\,v\ &,& (X. v)(dv)=X(dv)\,v\ .
\end{eqnarray} 
A covariant derivative $\nabla$ on $\Omega^{1}\Bbb T^2_q$ is specified by
$\nabla(du),\,\nabla(dv)\in \Omega^{1}\Bbb T^2_q\tens_{\Bbb T^2_q}
\Omega^{1}\Bbb T^2_q$.  We calculate the corresponding braiding from \ref{ppll} by
\begin{eqnarray}
\sigma(du\tens du) &=&\nabla(du.u)-\nabla(du).u \,=\,
\nabla(u.du)-\nabla(du).u \cr
&=& du\tens du+ u.\nabla(du)-\nabla(du).u\ ,\cr
\sigma(dv\tens du) &=&\nabla(dv.u)-\nabla(dv).u \,=\,q^{-1}\,
\nabla(u.dv)-\nabla(dv).u \cr &=& q^{-1}\, du\tens dv+ q^{-1}\,
u.\nabla(dv)-\nabla(dv).u\ ,\cr
\sigma(dv\tens dv) &=& dv\tens dv+ v.\nabla(dv)-\nabla(dv).v\
,\cr
\sigma(du\tens dv) &=& q\, dv\tens du+ q\,
v.\nabla(du)-\nabla(du).v\ .
\end{eqnarray}
To have a compatible interior product (see \ref{compatform}),
from \ref{lrcomp} we must have 
 $\Theta^{2}\Bbb T^2_q=\ker\wedge:\Omega^{1}\Bbb T^2_q\tens_{\Bbb T^2_q}
 \Omega^{1}\Bbb T^2_q\to \Omega^{2}\Bbb T^2_q$ contained in the +1 eigenspace
 of $\sigma$.  As $du\tens du$, $dv\tens dv$ and $du\tens dv+q\,dv\tens
 du$ are in $\Theta^{2}\Bbb T^2_q$, we deduce that
\begin{eqnarray}
    \nabla(du).u \,=\,u.\nabla(du) &,&
    \nabla(dv).v \,=\, v.\nabla(dv)\ , \cr
    q\,\nabla(dv).u-u.\nabla(dv) &=&
    q\,v.\nabla(du)-\nabla(du).v\ .
\end{eqnarray}
The general solution to this is, where the coefficients $r_{**}$
and $s_{**}$ are numbers,
\begin{eqnarray}\label{comcov}
    \nabla(du) &=& r_{uu}\, du\tens du.u^{-1}\,+\,
    r_{vu}\, dv\tens du.v^{-1}\,+\, r_{uv}\, du\tens dv.v^{-1}\,+\, r_{vv}\,
    dv\tens dv.u\,v^{-2}\ ,\cr
     \nabla(dv) &=& s_{vv}\, dv\tens dv.v^{-1}\,+\, s_{vu}\, dv\tens
     du.u^{-1}\,+\, s_{uv}\, du\tens dv.u^{-1}\,+\, s_{uu}\, du\tens
     du.v\,u^{-2}\ .
\end{eqnarray}
Putting this back into the equations for the braiding, we find
\begin{eqnarray}\label{tuu}
\sigma(du\tens du) 
&=& du\tens du\ ,\quad \sigma(dv\tens du)
 \,=\, q^{-1}\, du\tens dv\ ,\cr \sigma(dv\tens dv) &=& dv\tens
 dv\ ,\quad \sigma(du\tens dv) \,=\, q\, dv\tens du\ .
\end{eqnarray}
From (\ref{brad-eval}) and \ref{dubrad} we have
\begin{eqnarray}\label{etrr}
    \sigma(\partial_u\tens du) &=& du\tens \partial_u\
,\quad \sigma(\partial_v\tens du) \,=\, q\, du\tens \partial_v\ ,\cr
\sigma(\partial_v\tens dv) &=& dv\tens \partial_v\ ,\quad
\sigma(\partial_u\tens dv) \,=\, q^{-1}\, dv\tens \partial_u\ , \cr
\sigma^{-1}(\partial_u\tens \partial_u) &=& \partial_u\tens \partial_u
\ ,\quad \sigma^{-1}(\partial_v\tens \partial_u) \,=\, q^{-1}\,
\partial_u\tens \partial_v\ ,\cr \sigma^{-1}(\partial_v\tens
\partial_v) &=& \partial_v\tens \partial_v\ ,\quad
\sigma^{-1}(\partial_u\tens \partial_v) \,=\, q\, \partial_v \tens
\partial_u \ .
\end{eqnarray}
The paired left covariant derivative on $\vec \Bbb T^2_q$ can be calculated
using \ref{camab} as
\begin{eqnarray}
    \nabla(\partial_{u}) &=& -\,r_{uu}\, du.u^{-1}\tens \partial_{u}\,-\, q^{-1}r_{vu}\,
    dv.v^{-1}\tens \partial_{u}\,-\, q^{-1}s_{vu}\, dv.u^{-1}\tens\partial_{v}\,-\,
    s_{uu}\, du.v\,u^{-2}\tens \partial_{v}\ ,\cr 
    \nabla(\partial_{v}) &=& -\,s_{vv}\, dv.v^{-1}\tens \partial_{v}\,-\,
    q\,s_{uv}\, du.u^{-1}\tens \partial_{v}\,-\, q\,r_{uv}\,
    du.v^{-1}\tens \partial_{u}\,-\, r_{vv}\, dv.u\,v^{-2}\tens \partial_{u}\ .
\end{eqnarray}
From \ref{intdef}, the interior product on $\Omega^2 \Bbb T^2_q$ is
\[
\partial_u \lrcorner\,(du\wedge dv)\,=\, dv\ ,\quad
\partial_v \lrcorner\,(du\wedge dv)\,=\, -q\,du\ .
\]
From \ref{antidef} we see that $\partial_v
\tens\partial_u-q^{-1}\,\partial_u\tens\partial_v\in A^2 \Bbb T^2_q$, and
from \ref{deflie}, $\phi(\partial_v
\tens\partial_u-q^{-1}\,\partial_u\tens\partial_v)=0$.  The given
covariant derivatives restrict to $\Theta^{2}\Bbb T^2_q$, and the
covariant derivative on $\Omega^2 \Bbb T^2_q$ is
\begin{eqnarray}
\nabla(du\wedge dv) &=& du\tens(r_{uu}q^{-1}  + s_{uv})du\wedge dv.u^{-1}
\,+\, dv\tens(r_{vu}+q\,s_{vv})du\wedge dv.v^{-1}\ .
\end{eqnarray}
The Kroneker delta (see \ref{kroneker}) is $\delta=
du\tens\partial_{u}+dv\tens\partial_{v}$, and from (\ref{etrr}) $\dim \Bbb T^2_q=2$.

\begin{remark}
In summary; from the given differential calculus (\ref{nctrel})
we have derived a unique compatible braiding (\ref{tuu}) and a small
number of compatible covariant derivatives (\ref{comcov}). 
At this point, the reader may express some alarm: This is far too restrictive! 
Should we not be able to have a wide variety of covariant derivatives
on a given algebra?

Of course, the reader would be correct, we cannot restrict ourselves to such a small class
of possible covariant derivatives in the geometry of the noncommutative torus. 
However the problem is that we fixed the differential calculus (\ref{nctrel}).
We can indeed have many more covariant derivatives, but only if we allow different
differential calculi on our algebra. 
\end{remark}

\begin{remark}
In \cite{semiclass} the classical symplectic connection
associated to the semiclassical ($q\cong 1$) noncommutative torus and the given differential calculus
is computed. The braiding calculated from the connection
is shown to agree with the exact result here to lowest order in $q-1$. 
\end{remark}

\section{A noncommutative sphere}
    Following \cite{CHZ}, we describe a differential calculus on a
deformed sphere $S^2_q$ using a stereographic projection. The algebra on the 
coordinate chart of the projection is generated by $z$ and $\bar z$
with commutation relation $z\bar z=q^{-2}\bar zz+q^{-2}-1$. This can be made into a $C^*$
algebra with the involution $z^*=\bar z$. 
There is a left covariant (with respect to the action 
of $q$-deformed $SU_2$) differential calculus given by
\begin{eqnarray} \label{ncsph}
dz\wedge d\bar z=-q^{-2}\, d\bar z\wedge dz\ &,& z. dz=q^{-2}\, dz. z\ ,\quad 
z. d\bar z=q^{-2}\, d\bar z. z\ ,\cr  dz\wedge dz=d\bar z\wedge d\bar z=0  &,&
\bar z.  dz=q^{2}\, dz. \bar z\ ,\quad \bar z.  d\bar z=q^{2}\, d\bar
z.  \bar z\ .
\end{eqnarray} 
It will be convenient to set $R=\bar z\,z$, and 
as $R$ is a positive element in $S^2_q$, $R+1$ is invertible.  The
commutation relation can then be written $z(R+1)=q^{-2}(R+1)z$. 
For a function $f(R+1)$ we have $z\,f(R+1)=f(q^{-2}(R+1))\,z$.  Also we write
$z^1=z$ and $z^2=\bar z$.

\begin{propos}
The kernel of $\wedge:\Omega^1 S^2_q\tens_{S^2_q}\Omega^1 S^2_q\to \Omega^2 S^2_q$
is contained in the $+1$ eigenspace of $\sigma$ if
the covariant derivative is of the form
\begin{eqnarray*}
\nabla(dz) &=& q^{-2}\, (R+1)^{-1}\,\bar z. dz\tens dz\,+\,
(R+1)^{-1}\sum \big(g_{ij0}+q^6\,h_{ij1}\,z\big)\, dz^i\wedge dz^j\ ,
\cr \nabla(d\bar z) &=& q^4\,(R+1)^{-1}\,z. d\bar z\tens d\bar z\,+\,
(R+1)^{-1}\sum \big(h_{ij0}+h_{ij1}\,\bar z\big)\, dz^i\wedge dz^j\ ,
\end{eqnarray*}
where
$h_{ij0}$, $h_{ij1}$ and $g_{ij0}$ (for $1\le i,j\le 2$) are
constants. 
The braiding is 
\begin{eqnarray} \label{sphbra}
\sigma(dz\tens d z) &=& dz\tens d z \ ,\quad \sigma(d\bar z\tens d\bar z) 
\,=\, d\bar z\tens d\bar z\ ,\cr
\sigma(dz\tens d\bar z) &=& q^{-2}\, d\bar z\tens dz +(q^2-1) \sum h_{ij1}
\, dz^i\tens dz^j\ ,\cr
\sigma(q^{-2}\, d\bar z\tens dz) &=& dz\tens d\bar z -(q^2-1) \sum h_{ij1}
\, dz^i\tens dz^j\ .
\end{eqnarray}
\end{propos}
\proof
Begin with
\begin{eqnarray} 
dz\tens dz &=& \sigma(dz\tens dz)\,=\,\nabla(dz.z)-\nabla(dz).z \cr
&=& q^2\, \nabla(z.dz)-\nabla(dz).z\,=\, q^2\, dz\tens dz+q^2\,
z.\nabla(dz)-\nabla(dz).z\ .
\end{eqnarray} 
If we set $\nabla(dz)=\sum c_{ij}\, dz^i\wedge dz^j$, we have the
equations $z\,c_{ij}=q^2\,c_{ij}\,z$ for $(i,j)\neq(1,1)$, and
$1-q^2=q^2(z\,c_{11}-q^2\,c_{11}\,z)$.  The equation
$z\,c_{ij}=q^2\,c_{ij}\,z$ has solution $c_{ij}=(R+1)^{-1}\,g_{ij}(z)$
for any non-singular function $g_{ij}$.  Now the equation for $c_{11}$
has solution
\[
c_{11}\,=\,q^{-2}\,z^{-1}\,+\,(R+1)^{-1}\,f_{11}(z)\ .
\]
Unfortunately $z$ is not invertible, but we can write
$(1-(R+1)^{-1})\,z^{-1}=(R+1)^{-1}\,\bar z$, and then
\begin{eqnarray}
c_{11} &=&  q^{-2}\, (R+1)^{-1}\,\bar z\,+\, (R+1)^{-1}\,g_{11}(z)  \ ,
\end{eqnarray}
where $g_{11}(z)$ is non-singular. Next,
\begin{eqnarray} 
d\bar z\tens d\bar z &=& \sigma(d\bar z\tens d\bar z)\,=\,\nabla(d\bar z.\bar z)
-\nabla(d\bar z).\bar z \cr
&=& q^{-2}\, \nabla(\bar z.d\bar z)-\nabla(d\bar z).\bar z\,=\, q^{-2}
\, d\bar z\tens d\bar z+q^{-2}\, \bar z.\nabla(d\bar z)-\nabla(d\bar z).\bar z\ .
\end{eqnarray} 
We set $\nabla(d\bar z)=\sum e_{ij}\, dz^i\wedge dz^j$, and get the equations
$\bar z\,e_{ij}=q^{-2}\,e_{ij}\,\bar z$ for $(i,j)\neq (2,2)$ and 
\[
1-q^{-2}\,=\,q^{-2}\,\bar z\,\,e_{22}-q^{-4}\,e_{22}\,\bar z\ .
\]
Now we use the result $\bar z\,(R+1)^{-1}=q^{-2}\,(R+1)^{-1}\,\bar z$
to find that $e_{ij}=(R+1)^{-1}\, h_{ij}(\bar z)$ for $(i,j)\neq (2,2)$ and 
\[
e_{22}\,=\,q^4\,(R+1)^{-1}\,z\,+\,(R+1)^{-1}\,h_{11}(\bar z)\ .
\]
Now we use the fact that $d z\tens d\bar z+ q^{-2}\,d\bar z\tens dz$ is 
in the kernel of $\wedge$:
\begin{eqnarray} \label{wedgesum}
    \sigma(d z\tens d\bar z) &=& \nabla(d z.\bar
z)-\nabla(d z).\bar z \,=\, q^{-2}\,\nabla(\bar z.dz)-\nabla(d z).\bar
z \cr &=& q^{-2}\,d\bar z\tens dz+q^{-2}\,\bar z.\nabla(dz)-\nabla(d
z).\bar z\ ,\cr q^{-2}\, \sigma(d \bar z\tens dz) &=&
q^{-2}\,\nabla(d\bar z.z)-q^{-2}\, \nabla(d\bar z).z \,=\,
\nabla(z.d\bar z)-q^{-2}\,\nabla(d\bar z).z \cr &=& dz\tens d\bar
z+z.\nabla(d\bar z)-q^{-2}\,\nabla(d\bar z).z \ .
\end{eqnarray} 
Then, using $z(R+1)^{-1}=q^{2}(R+1)^{-1}z$,
\begin{eqnarray} 
q^{-2}\,\bar z.\nabla(dz)-\nabla(d z).\bar z &=& \sum (q^{-2}\bar z\
,c_{ij}-q^{-4}\,c_{ij}\bar z)
dz^i\tens dz^j \cr
&=&q^{-4}\,(R+1)^{-1}\sum (\bar z\,g_{ij}(z)-g_{ij}(z)\bar z)
dz^i\tens dz^j\ ,\cr
q^{-2}\,\nabla(d\bar  z).z-z.\nabla(d\bar  z) &=& \sum(q^2\,e_{ij}
\,z-z\,e_{ij})dz^i\tens dz^j\cr
&=& q^2\,(R+1)^{-1}\sum (h_{ij}(\bar z)\,z-z\,h_{ij}(\bar z))dz^i\tens dz^j\ ,
\end{eqnarray} 
and we find that
\[
q^6(h_{ij}(\bar z)\,z-z\,h_{ij}(\bar z))\,=\,\bar z\,g_{ij}(z)-g_{ij}(z)\bar z\ .
\]
A solution to this is given by constants $h_{ij0}$, $h_{ij1}$ and
$g_{ij0}$ (for $1\le i,j\le 2$), when
\[
h_{ij}(\bar z)\,=\, h_{ij0}+h_{ij1}\,\bar z\ ,\quad
g_{ij}(z)\,=\, g_{ij0}+q^6\,h_{ij1}\,z\ .
\]
The braiding is calculated from (\ref{wedgesum}).\quad$\square$

\begin{propos} \label{list} Assume that $q$ is nonzero and $q^4\neq 1$. 
The cases for which the braiding (\ref{sphbra}) gives an interior
product which is compatible with the differential calculus
(\ref{compatform}) are

a)\quad $h_{221}=h_{211}=0$, $h_{121}=1/(q^{2}-1)$.

b)\quad $h_{111}=h_{121}=0$, $h_{211}=1/(q^{2}-q^{4})$.

c)\quad $h_{111}=h_{221}=0$, $h_{211}\,h_{121}=0$.

d)\quad $h_{111}=h_{221}=0$,  $h_{121}=1/(q^{2}-1)$, $h_{211}=1/(q^{2}-q^{4})$.

\end{propos}
\proof The only case left to check is that $T_3 \Theta^3 S^2_q\subset \Omega^1 S^2_q \tens_{S^2_q}
\Theta^2 S^2_q$, where 
$\Theta^3 S^2_q$ is all of $\Omega^1 S^2_q \tens_{S^2_q}\Omega^1 S^2_q \tens_{S^2_q}\Omega^1 S^2_q$,
 and explicit calculation gives the answer.\quad$\square$

\begin{propos} Assume that $q$ is nonzero and $q^4\neq 1$. 
The cases for which the generalised braiding (\ref{sphbra}) actually satisfies
the braid relation are

a)\quad $h_{111}=h_{221}=0$, $h_{211}\,h_{121}=0$.

b)\quad $h_{111}=h_{221}=0$,  $h_{121}=1/(q^{2}-1)$, $h_{211}=1/(q^{2}-q^{4})$.

\end{propos}

\begin{propos} Assume that $q$ is nonzero and $q^4\neq 1$. 
The condition for the braiding (\ref{sphbra}) to be invertible is that 
$(q^{2}-1)(h_{121}-h_{211}\,q^{2})\neq 1$, and the cases for which $\sigma^2$ is the
identity are

a)\quad $h_{121}-h_{211}\,q^{2}=0$.

b)\quad $h_{111}=h_{221}=0$,  $h_{121}=1/(q^{2}-1)$, $h_{211}=1/(q^{2}-q^{4})$.

\end{propos}

\begin{propos} The vector fields are generated by $\partial_z$ and $\partial_{\bar z}$,
where $\partial_z(dz)=\partial_{\bar z}(d\bar z)=1$ and 
$\partial_z(d\bar z)=\partial_{\bar z}(dz)=0$. Then the `Lie bracket'
is $\phi(\partial_{\bar z}\tens \partial_{ z}-q^2\,\partial_{z}\tens \partial_{\bar  z})=0$,
and the braiding
$\sigma:\vec\tens\Omega^1 \to \Omega^1\tens \vec$ is given by
\begin{eqnarray} 
\sigma(\partial_{z}\tens d z) &=& dz\tens \partial_{z} +
\frac {h_{111}\,(1-q^{2})}x dz\tens \partial_{\bar z} + 
\frac{h_{211}\,(q^{2}-q^{4})}{x} d\bar z\tens \partial_{\bar z}\ ,\cr
\sigma(\partial_{z}\tens d \bar z) &=& \frac{h_{111}\,(q^{4}-q^{2})}{x}dz\tens
\partial_{z} + \frac{h_{121}\,(q^{4}-q^{2})-q^{2}}{x} d\bar z\tens \partial_{z}\ ,\cr
\sigma(\partial_{\bar z}\tens d z) &=&
\frac{h_{211}\,(1-q^{2})-1/q^{2}}{x} dz\tens \partial_{\bar z} +
\frac{h_{221}\,(1-q^{2})}{x} d\bar z\tens \partial_{\bar z}\ ,\cr
\sigma(\partial_{\bar z}\tens d \bar z) &=&\frac{h_{121}\,(q^{2}-1)}{x} dz\tens
\partial_{z} \,+\,\frac{h_{221}\,(q^{4}-q^{2})}{x} d\bar z\tens \partial_{z} +
d\bar z\tens \partial_{\bar z}\ ,
\end{eqnarray}
where $x=(q^{2}-1)(h_{121}-q^{2}\,h_{211})-1$.  The differential dimension is
\[
\dim S^2_q =
\frac{x(x-1)}{x^{2}+q^{2}\,(q^{2}-1)^{2}\,(h_{121}\,h_{211}-h_{111}\,h_{221})}\ .
\]
\end{propos}

\begin{remark}
In summary; we now have a wider variety of possibilities.  The differenital
calculus does not uniquely specify the braiding.  There are
generalised braidings (\ref{list} (a) with $h_{111}\neq 0$ and (b)
with $h_{221}\neq 0$) which are compatible with the differenital
calculus but do not satisfy the braid relation.  There are generalised
braidings (\ref{list} (c) with exactly one of $h_{121}$ or $h_{211}$
vanishing) which are compatible with the differenital calculus and
satisfy the braid relation, but do not square to the identity
and give (in general) fractional differential dimension.
\end{remark}

\section{Curvature and Torsion}
Using vector fields, we can define the curvature and torsion by
some remarkably familiar classical formulae, rather than the usual
noncommutative formulae using forms.  Remember that $\phi$
defined in \ref{deflie} is the analogue of the Lie bracket.

\begin{defin} Given a left $M$-covariant derivative $\nabla$ on a left $M$-module 
    $E$, define the curvature
$R:A^2M\tens E\to E$ as
\[
R(X\tens Y)(e)\,=\,\nabla_X\,\nabla_Y\,(e)\,-\,\nabla_{\phi(X\tens Y)}\,(e)\ .
\]
\end{defin}

\begin{propos} The curvature descends to a well defined left $M$-module map
    $R:\pi A^2M\tens_M E\to E$, where $\pi:\vec M\tens \vec M\to \vec
    M\tens_M \vec M$ is the quotient map.
\end{propos}
\proof The left module property is quite simple. Next, for $X\tens Y\in A^2M$:
\begin{eqnarray*}
R(X\tens m.Y)(e) &=& \nabla_X \nabla_{m.Y}(e)\,-\,\nabla_{\phi(X\tens m.Y)}(e) \cr
&=& \nabla_X (m.\nabla_{Y}(e))\,-\,\nabla_{\phi(X.m\tens Y)}(e)\,-\,\nabla_{D_X(m).Y}(e) \cr
&=& R(X.m\tens Y)(e)\ .
\end{eqnarray*}
Next, using \ref{prop17},
\begin{eqnarray*}
R(X\tens Y)(m.e) &=& \nabla_X \nabla_{Y}(m.e)\,-\,\nabla_{\phi(X\tens Y)}(m.e) \cr
&=& \nabla_X(D_Y(m).e+\nabla_{Y.m}(e))\,-\,\phi(X\tens Y)(dm).e\,-\,
\nabla_{\phi(X\tens Y).m}(e) \cr
&=& R(X\tens Y.m)(e) \,+\, \nabla_X(D_Y(m).e) \,-\, \phi(X\tens Y)(dm).e\,-\,
\nabla_{X.D_Y(m)}(e)\cr
&=& R(X\tens Y.m)(e)\ .\quad\square
\end{eqnarray*}

\begin{defin} Given a bimodule connection on $\Omega^{1}M$, we define the torsion
 $T:A^2M\to {\rm Vec}\,M$ as
$
T(X\tens Y)\,=\,\nabla_X(Y)\,-\,\phi(X\tens Y)
$.
\end{defin}

\begin{propos} The torsion descends to a well defined left module map
    $T:\pi A^2M\to {\rm Vec}\,M$.
\end{propos}
\proof Using \ref{prop17},
\begin{eqnarray*}
T(X.m\tens Y) &=& \nabla_{X.m}(Y)\,-\,\phi(X.m\tens Y) \cr
&=& \nabla_{X.m}(Y)\,-\,\phi(X\tens m.Y)\,+\, D_X(m).Y\,=\,T(X\tens m.Y)\ .
\end{eqnarray*}
For the left module map condition,
\begin{eqnarray*}
T(m.X\tens Y) &=& \nabla_{m.X}(Y)\,-\,\phi(m.X\tens Y) \,=\,
m.T(X\tens Y)\ .\quad\square
\end{eqnarray*}

\section{Classical exponentiation and parallel transport}

We consider the usual point based definition of differential geometry
on a manifold $M$, and translate 
it into a form more amenable to non-commutative geometry.
The directional derivative notation is used,
where $f'(x;v)$ is the derivative of the function at $x$ along the
vector $v$.  We begin with a result about exponentiating a known time
dependent vector field into a diffeomorphism:

\begin{propos} \label{expa} Given a vector field $w(t)\in{\rm Vec}\,M$ for
$t\in\Bbb R$, define a function $p:M\times \Bbb R\to M$ by $p(x,0)=x$
and $\dot p(x,t)\,=\,w(t)(p(x,t))$.  Now define functions $J,K:\Bbb
R\to L(C^\infty(M),C^\infty(M))$ (the set of linear maps from
$C^\infty(M)$ to itself) by
\[
J(t)(f)(x) \,=\, f'(x;w(x,t))\ ,\quad 
K(t)(f)(x)\,=\,f(p(x,t))\ .
\] 
Then $K(t)\in  L(C^\infty(M),C^\infty(M))$ is the solution to the differential equation
\[
\dot K(t)\,=\, K(t)\circ J(t)\ ,\quad K(0)={\rm id}\ .
\]
\end{propos}
\proof By differentiating the definition of $K(t)$ with respect to $t$,
\[
\dot K(t)(f)(x)\,=\,f'(p(x,t);w(p(x,t),t))\,=\,
K(t)(J(t)(f))(x)
\ .\quad\square
\]

Now we consider parallel transport in a bundle $E$ with connection $\nabla$, along a curve 
which is given by exponentiation of a time dependent vector field. 
The connection is specified in our coordinate system by Christoffel
symbols $\Gamma$.

\begin{propos}
Take a curve $p(x,t)$ starting at $x$ given by exponentiating the time dependent 
vector field $w(t)$. Along the curve take
$s(x,t)\in E_{p(x,t)}$ which is a solution to the parallel transport equation
$\dot s(x,t)\,+\,\Gamma(p(x,t);\dot p(x,t),s(x,t))\,=\,0$. Define
a time dependent section $c(t)$ of $E$ by $c(t)(p(x,t))\,=\,s(x,t)$. 
Then $c(t)$ obeys the first order differential equation $\dot
c(t)\,=\,-\,\nabla_{w(t)}c(t)$.
\end{propos}
\proof
Differentiating the definition of $c$ with respect to $t$ we find
\[
\dot s(x,t)\,=\,\dot c(t)(p(x,t))\,+\,c(t)'(p(x,t);\dot p(x,t))\,=\,-\,
\Gamma(p(x,t);\dot p(x,t),s(x,t))\ .
\]
We can rerarrange this to give 
\[
\dot c(t)(p(x,t))\,=\,-\,\nabla_{\dot p(x,t)}c(t)(p(x,t))\ .\quad\square
\]

\medskip 
Now we can use a connection on
the tangent bundle $TM$ to define geodesics on $M$.

\begin{cor}
Suppose that the vector field $c(t)\in {\rm Vec}\,M$ is parallel transported
along curves which are given by exponentiating $c(t)$ itself. Then $c(t)$ obeys the 
first order non-linear differential equation $\dot
c(t)\,=\,-\,\nabla_{c(t)}c(t)$.
\end{cor}

\medskip The reader will notice that the geodesic equation is, unlike the 
parallel transport equation, non-linear. First order linear equations tend (sweeping
much under the carpet) to have solutions which can be extended for all
time, wheras non-linear equations can easily have solutions which blow
up at finite time.  This phenomenon is well known in classical
geometry, in fact a manifold is called complete just when its
geodesics can be extended for all time.

\section{Non-commutative vector fields and parallel transport}
Now we translate the ideas of the last section into the non-commutative regime.
There is a problem with exponentiation, if $X(t)$ is not a derivation then its exponentiation
is not an algebra map. A partial answer is given in this section.
We denote by $L(A,B)$ the linear maps from $A$ to $B$.

\begin{defin} A time dependent vector field $X(t)\in {\rm Vec}\,(M)$ exponentiates
to give $K^n_X(t)\in L(\Omega^n M,\Omega^n M)$
(for all $n\ge 0$) defined by
\[
\dot K^n_X(t)\,=\, K^n_X(t)\circ {\cal L}_{X(t)}\ ,\quad K^n_X(0)={\rm id}\ .
\]
For $n=0$ in the classical case, this gives the same result as \ref{expa}.
\end{defin}

\begin{propos}\label{coch4} The sequence of maps $K^*_X(t):\Omega^* M\to
\Omega^* M$ is a cochain complex map, i.e.\ $d\circ
K^n_X(t)=K^{n+1}_X(t)\circ d:\Omega^n M\to \Omega^{n+1} M$.
\end{propos}
\proof This is given by the uniqueness of solutions to first order 
 equations.  Using \ref{lieprops},
\begin{eqnarray*}
(\dot K^{n+1}_X(t)\circ d)(\xi) &=& (K^{n+1}_X(t)\circ{\cal L}_{X(t)})(d\xi)\,=\,
(K^{n+1}_X(t)\circ d)({\cal L}_{X(t)}(\xi))\ ,\cr
(d\circ \dot K^n_X(t))(\xi) &=& (d\circ K^n_X(t))({\cal L}_{X(t)}(\xi))\ .
\end{eqnarray*}
Then $d\circ K^n_X(t)$ and $K^{n+1}_X(t)\circ d$ are both solutions to the
differential equation $\dot U(t)=U(t)\circ {\cal L}_{X(t)}$ with
initial condition $U(0)=d$ for $U(t):\Omega^n M\to \Omega^{n+1} M$.
\quad$\square$

\begin{propos} The maps $\dot K_X^n(t)$ are cochain homotopic to 0 via the 
    cochain homotopy
$h_X^n(t)=K_X^n(t)\circ (X(t)\lrcorner):\Omega^{n+1}M\to\Omega^n M$. 
\end{propos}
\proof We use the definition of the Lie derivative \ref{liederdef}
and \ref{coch4} to write
\[
\dot K_X^{n+1}(t) \,=\, K_X^{n+1}(t)\circ(d\circ (X(t)\lrcorner)+(X(t)\lrcorner)\circ d)
\,=\, d\circ(K_X^n\circ (X(t)\lrcorner))+(K_X^{n+1}\circ (X(t)\lrcorner))\circ d
\ .\quad\square
\]

\begin{defin} Let $E$ be a left $M$-module with connection $\nabla$, and
take a time dependent vector field $X(t)\in {\rm Vec}\,(M)$.  Then
$c(t)\in E$ is parallel transported along the exponentiation of $X(t)$
if it obeys the first order differential equation $\dot
c(t)\,=\,-\,\nabla_{X(t)}c(t)$.
\end{defin}

\begin{defin} Given a connection $\nabla$ on ${\rm Vec}\,M$, $c(t)\in {\rm Vec}\,M$ 
    is parallel transported along the 
exponentiation of $c(t)$ if it 
obeys the first order differential equation $\dot c(t)\,=\,-\,\nabla_{c(t)}c(t)$. 
\end{defin}

\section{Exponentials of vector fields on the non-commutative torus}
We will compare the exponentials of the time independent vector fields
$u.\partial_u$ (which is a bimodule map) and $\partial_u$ (which is
only a right module map) on $\Bbb T^{2}_{q}$.

\begin{lemma} We have $d(v^r u^s)=r\, dv.v^{r-1}u^s + s\,q^{-r}\, du.v^r u^{s-1}$. 
    This then gives
\[
{\cal L}_{\partial_u}(v^r u^s) \,=\, s\,q^{-r}\, v^r u^{s-1}\ ,\quad
{\cal L}_{u\partial_u}(v^r u^s) \,=\, s\, v^r u^{s}\ .
\]
\end{lemma}

\begin{propos} On $\Omega^0 \Bbb T^2_q$ we get
\[
\exp(t\,{\cal L}_{\partial_u})(v^r u^s) \,=\, (1+t\,u^{-1})^s\,v^r\, u^s\ ,\quad
\exp(t\,{\cal L}_{u\partial_u})(v^r u^s) \,=\, v^r\, u^s\,e^{st}\ .
\]
\end{propos}
\proof For the difficult case, first iterate the Lie derivative to get
\begin{eqnarray*}
    ({\cal L}_{\partial_u})^n(v^r u^s) &=&
s(s-1)\dots(s-n+1)\,q^{-nr}\,v^r\,u^{s-n} \cr
 &=&
s(s-1)\dots(s-n+1)\,u^{-n}\,v^r\,u^{s}
\ ,
\end{eqnarray*}
and then use the binomial expansion
\[
\sum_{n\ge 0}\frac{s(s-1)\dots(s-n+1)}{n!}\,(t\,u^{-1})^n\,=\,
(1+t\,u^{-1})^s\ .\quad\square
\]

\begin{lemma} On $\Omega^1 \Bbb T^2_q$ we get
\begin{eqnarray*}
\exp(t\,{\cal L}_{\partial_u})(du.v^r u^s+dv.v^n u^m) &=&
(1+t\,u^{-1})^s.du.v^r u^s\,\,
+\,(1+t\,\,u^{-1})^m.dv.v^n u^m\ ,\cr
\exp(t\,{\cal L}_{u\partial_u})(du.v^r u^s+dv.v^n u^m) &=& 
e^{(s+1)t}\,du.v^r u^s\,+\, e^{mt}\,dv.v^n u^m\ .
\end{eqnarray*}
\end{lemma}
\proof From the following equations:
\begin{eqnarray*}
{\cal L}_{\partial_u}(du.v^r u^s+dv.v^n u^m) &=& 
d(\partial_u\lrcorner(du.v^r u^s+dv.v^n u^m)) \cr
&&-\,\partial_u\lrcorner(du\wedge d(v^r u^s) + dv\wedge d(v^n u^m)) \cr
&=& d(v^r u^s) \,-\,\partial_u\lrcorner(du\wedge dv).r\,v^{r-1} u^s\,-\,
\partial_u\lrcorner(dv\wedge du).m\,q^{-n}\,v^nu^{m-1} \cr
&=&  s\,q^{-r}\, du.v^r u^{s-1}\,+\, dv.m\,q^{-n-1}\,v^nu^{m-1}\ ,\cr
{\cal L}_{u\partial_u}(du.v^r u^s+dv.v^n u^m) &=& 
d(u\,\partial_u\lrcorner(du.v^r u^s+dv.v^n u^m)) \cr
&&-\,u\,\partial_u\lrcorner(du\wedge d(v^r u^s) + dv\wedge d(v^n u^m)) \cr
&=& q^r\,d(v^r u^{s+1}) \,-\,u\,\partial_u\lrcorner(du\wedge
dv).r\,v^{r-1} u^s\cr && -\, u\,\partial_u\lrcorner(dv\wedge
du).m\,q^{-n}\,v^nu^{m-1} \cr &=& (s+1)\, du.v^r u^{s}\,+\,
m\,dv.v^nu^{m}\ .\quad\square
\end{eqnarray*}

\begin{propos} On $\Omega^2 \Bbb T^2_q$ we get
\begin{eqnarray*}
\exp(t\,{\cal L}_{\partial_u})(du\wedge dv.v^r u^s) 
&=& (1+t\,u^{-1})^s.du\wedge dv.v^r u^s\ ,\cr
\exp(t\,{\cal L}_{u\partial_u})(du\wedge dv.v^r u^s) 
&=& du\wedge dv.v^r u^s\,e^{(s+1)t}\ .
\end{eqnarray*}
\end{propos}
\proof From the following equations:
\begin{eqnarray*}
{\cal L}_{\partial_u}(du\wedge dv.v^r u^s) &=&
d(\partial_u\lrcorner(du\wedge dv).v^r u^s) \cr
&=& d(dv.v^r u^s)\,=\,-\,dv\wedge d(v^r u^s) \cr
&=& -\,dv\wedge du.s\,q^{-r}\,v^r u^{s-1}\,=\,
du\wedge dv.s\,q^{-r-1}\,v^r u^{s-1}\ ,\cr
{\cal L}_{u\partial_u}(du\wedge dv.v^r u^s) &=&
d(u.\partial_u\lrcorner(du\wedge dv).v^r u^s) \cr
&=& d(u.dv.v^r u^s)\,=\,q^{r+1}\,d(dv.v^r u^{s+1})\,=\,
-\,q^{r+1}\,dv\wedge d(v^r u^{s+1}) \cr
&=& -\,q\,dv\wedge du.(s+1)\,v^r u^s\,=\,
(s+1)\,du\wedge dv.v^r u^s\ .\quad\square
\end{eqnarray*}

\section{Exponentiation and Hopf algebra coactions}
It would be somewhat premature for me to claim that these 
exponentials of Lie derivatives really
 were significant in the non-commutative context, just 
 because they reduce to the correct construction in
the commutative case. Thus I would like to present some 
non-commutative supporting evidence.

Given a differentiable action of a Lie group on a manifold, an element
of the Lie algebra gives a vector field on the manifold. 
Exponentiation of this vector field gives a diffeomorphism which is
just action by the exponential of the Lie algebra element as an
element of the Lie group.  In this section I show an analagous result for
Hopf algebra coactions on algebras.  I shall use the Sweedler notation
$\Delta(h)=h_{(1)}\tens h_{(2)}$ for coproducts.

\begin{defin} Suppose that a Hopf algebra $H$ is given a 
differentiable structure
so that the coproduct $\Delta:H\to H\tens H$ is differentiable, where
$H\tens H$ is given the tensor product differential structure (see \ref{tensdef}). 
The braided Lie algebra of $H$ \cite{brLie} is defined as
\[
\gh=\{\alpha:\Omega^1 H\to k:\alpha(\xi.h)=\alpha(\xi).\eps(h)
\quad \forall h\in H\ \forall\xi\in \Omega^{1}M\}\ .
\]
\end{defin}

\begin{remark} This idea of differentiability is really the 
same as the more usual
idea of bicovariance of the differential calculus. Given the 
existence of $\Delta_*$,
we define right and left coactions of $H$ on $\Omega^1 H$ by 
$\rho=\Pi_{1}\Delta_{*}$
and $\mu=\Pi_{2}\Delta_{*}$ respectively. The fact that these are coactions can be 
checked from the tensorial property in \ref{tensdef}. 
\end{remark}

\begin{remark} See \cite{worondiff} for more on the differential calculus on
Hopf algebras. Note that the condition that $\Delta:H\to H\tens H$ is differentiable
is the same as requiring the bicovariance of the calculus, 
where the right and left coactions are $\Pi_1\circ \Delta_*:\Omega^1 H\to
\Omega^1 H\tens H$ and $\Pi_2\circ \Delta_*:\Omega^1 H\to
H\tens \Omega^1 H$.
\end{remark}

\begin{propos} There is a 1-1 correspondence between $\gh$ and left
$H$-covariant vector fields on $H$ given by $\alpha\in\gh$ mapping to
$L_\alpha=(\id\tens\alpha)\Pi_2 \Delta_*:\Omega^1 H\to H$, and a
vector field $X$ mapping to $\eps\circ X\in\gh$.
\end{propos}
\proof For a vector field $X$,
$\eps(X(\xi.h))=\eps(X(\xi).h)=\eps(X(\xi))\,\eps(h)$, so $\eps\circ
X\in\gh$.  Also
\[
(\id\tens\alpha)\Pi_2
\Delta_*(dh.a)\,=\, h_{(1)}\,a_{(1)}\,\alpha(dh_{(2)}.a_{(2)}) \,=\,
h_{(1)}\,a_{(1)}\,\eps(a_{(2)})\,\alpha(dh_{(2)})\,=\,h_{(1)}\,\alpha(dh_{(2)})
\,a\ ,
\]
so $L_{\alpha}$ is a right module map, i.e.\ a vector field on $H$. 
To check that $L_{\alpha}$ is left invariant,
\begin{eqnarray*}
    (\id\tens L_{\alpha})\mu(dh.a) &=& h_{(1)}a_{(1)}\tens L_{\alpha}
    (dh_{(2)}.a_{(2)})\,=\, h_{(1)}a_{(1)}\tens h_{(2)(1)}\,a_{(2)}\,
    \alpha(dh_{(2)(2)})\ ,\cr
    \Delta\circ L_{\alpha}(dh.a) &=& \alpha(dh_{(2)})\,\Delta(h_{(1)}\,a)
    \,=\,\alpha(dh_{(2)})\,h_{(1)(1)}\,a_{(1)}\tens h_{(1)(2)}\,a_{(2)}\ ,
\end{eqnarray*}
and these are the same by coassociativity.  To check the 1-1 correspondence,
\[
\eps(L_{\alpha}(dh)) \,=\, \eps(h_{(1)})\,\alpha(dh_{(2)})\,=\,\alpha(dh)\ ,
\]
and finally the more difficult bit.  From the discussion above
\[
L_{\eps\circ X}(dh)\,=\, h_{(1)}\,\eps(X(dh_{(2)})\ ,
\]
and if $X$ is left invariant then $\Delta(X(dh))=h_{(1)}\tens
X(dh_{(2)})$, so 
\[
L_{\eps\circ X}(dh)\,=\, (\id \tens\eps)\Delta(X(dh))\,=\,X(dh)\ .\quad\square
\]

\begin{propos}\label{expp} The exponentiation of the time independent
vector field $L_\alpha$ on $H$ is
\[
\exp(t{\cal L}_{L_\alpha})(h)\,=\,
\sum_{r\ge 0} \frac{t^r}{r!}\,
h_{(1)}\,\alpha(dh_{(2)})\dots\alpha(dh_{(r+1)})\ .
\]
\end{propos}
\proof First calculate the Lie derivative ${\cal L}_{L_\alpha}(h)=
{L_\alpha}(dh)=h_{(1)}\,\alpha(dh_{(2)})$.  Iterating this and using
coassociativity we find $({\cal
L}_{L_\alpha})^2(h)=h_{(1)}\,\alpha(dh_{(2)})\,\alpha(dh_{(3)})$ etc. 
\quad $\square$

\begin{defin} We define the exponential of $\alpha\in \gh$ to be
    an element $\exp(\alpha)$ of the dual $H^{*}$.  For $h\in H$,
    $\exp(\alpha)(h)$ is defined to be the counit $\eps$ applied to
    $\exp({\cal L}_{L_\alpha})(h)$, i.e.\ from \ref{expp}
    \[
    \exp(\alpha)(h)\,=\,\sum_{r\ge 0} \frac{1}{r!}\,
    \alpha(dh_{(1)})\dots\alpha(dh_{(r)})\ .
    \]
    As a brief check that this corresponds to the classical construction,
    the Lie algebra $\gg$ of a Lie group $G$ in our setting corresponds to 
    the Hopf algebra $k(G)$ of functions on $G$. The exponential of an
    element of $\gg$ is in the dual algebra, the group algebra $kG$.
\end{defin}

Now we turn to a differentiable left coaction $\lambda$ of a Hopf algebra
$H$ on an algebra $M$. We suppose that $M$ is a left $H$-comodule
algebra, i.e.\ $\lambda:M\to H\tens M$ is an algebra map (with the tensor product
algebra structure) and (if $M$ is unital) $\lambda(1_M)=1_H\tens 1_M$. 

\begin{propos} Given a
left $H$-comodule
algebra $M$ with differentiable left $H$-coaction $\lambda:M\to H\tens
M$, there is a map $\Lambda :\gh\to \vec M$ given by
$\Lambda(\alpha)= (\alpha\tens\id)\Pi_1\lambda_*$.
\end{propos}
\proof If we write the left coaction as $\lambda(m)=m_{[-1]}\tens
m_{[0]}$, then
\begin{eqnarray*}
\lambda_*	(dn.m) &=& d(\lambda(n)).\lambda(m)\,=\,
d(n_{[-1]}	\tens n_{[0]}	).(m_{[-1]}	\tens m_{[0]}	)\ ,\cr
\Pi_1	\lambda_*	(dn.m) &=& dn_{[-1]}.m_{[-1]}\tens ,n_{[0]}m_{[0]}\ ,\cr
\Lambda(\alpha)(dn.m) &=& \alpha(dn_{[-1]}.m_{[-1]})\,n_{[0]}m_{[0]}\,=\,
\alpha(dn_{[-1]})\,n_{[0]}\,\eps(m_{[-1]})\,m_{[0]}\,=\,
\Lambda(\alpha)(dn).m\ ,
\end{eqnarray*}
so $\Lambda(\alpha)$ is a right module map, i.e.\ a vector field on $M$.
\quad$\square$

\begin{propos} The exponential of the time independent vector field
$\Lambda(\alpha)$ on $M$ is
\[
\exp(t{\cal L}_{\Lambda(\alpha)})(m)\,=\, \sum_{r\ge 0} \frac{t^r}{r!}\,
m_{[0]}\,\alpha(dm_{[-1](1)})\dots\alpha(dm_{[-1](r)})\ .
\]
\end{propos}
\proof First calculate the Lie derivative 
${\cal
L}_{\Lambda(\alpha)}(m)={\Lambda(\alpha)}(dm)=\alpha(dm_{[-1]})\,m_{[0]}$.
 Iterating this and using the coaction property, $({\cal
L}_{\Lambda(\alpha)})^2(m)=\alpha(dm_{[-1](1)})\,
\alpha(dm_{[-1](2)})\,m_{[0]}$ etc.  \quad $\square$

\begin{theorem} We have the following relation between the exponential of
the vector field on $M$ generated by an element $\alpha\in\gh$ and the
exponential of $\alpha$ as an element of $H^{*}$:
\[
\exp({\cal L}_{\Lambda(\alpha)})\,=\,(\exp(\alpha)\tens\id)\circ
\lambda:M\to M\ .
\]
\end{theorem}
\proof Directly from the preceeding results.\quad $\square$

\end{document}